\documentclass[11pt,reqno, a4wide]{amsart}
\usepackage{amssymb}
\usepackage{amsmath}
\usepackage{amsthm}
\usepackage{mathptmx}
\usepackage{color}
\usepackage[pdftex]{graphicx}
\usepackage{hyperref}
\usepackage{a4wide}
\usepackage[latin1]{inputenc}

\newtheorem{teorema}{Theorem}

\newtheorem{defini}[teorema]{Definition}
\newtheorem{lema}[teorema]{Lemma}

\newtheorem{prop}[teorema]{Proposition}
\newtheorem{obs}[teorema]{Remark}

\newcommand{\R}{\mathbb{R}}%cria um novo comando de nome \R

\newcommand{\vol}{\mathrm{vol}}

\newcommand{\Spec}{\mathrm{Spec}}
\newcommand{\N}{\mathbb{N}}
\newcommand{\Ree}{\mathrm{Re}}
\newcommand{\be}{\begin{equation}}
\newcommand{\ee}{\end{equation}}
\newcommand{\bee}{\begin{equation*}}
\newcommand{\eee}{\end{equation*}}

\newcommand{\bH}{\mathbf{H}}

\newcommand{\whitebox}{$\Box$}
\newenvironment{prova}[1][Proof]{\textbf{#1.} }{\hfill\whitebox\bigskip}
\newenvironment{dem}[1][Proof]{\textbf{#1.} }{\hfill\whitebox\bigskip}
\numberwithin{equation}{section}

\newcommand{\supp}{\mathrm{supp}}

\begin{document}

\title[Quadratic  Schr\"odinger  system]{ Local and global well-posedness for  a quadratic  Schr\"odinger  system on  Zoll manifolds}
\author{ Marcelo Nogueira}
\address{Department of Mathematics, State University of Campinas, 13083-859, Campinas, SP, Brazil}
\email{marcelonogueira19@gmail.com }
\thanks{M. Nogueira was supported by CAPES and  CNPq, Brazil.}
\author{Mahendra Panthee}
\address{Department of Mathematics, State University of Campinas, 13083-859, Campinas, SP,  Brazil}
\email{mpanthee@ime.unicamp.br}
\thanks{M. Panthee was partially supported by CNPq (308131/2017-7) and FAPESP (2016/25864-6) Brazil.}


\begin{abstract}
We consider the initial value problem (IVP) associated to a quadratic  Schr\"odinger system
\begin{equation*}
\begin{cases}
i \partial_{t} v \pm \Delta_{g} v -  v = \epsilon_{1} u \bar{v}, & t \in \R,\; x \in M, \\ 
i \sigma \partial_{t} u \pm  \Delta_{g} u - \alpha u = \frac{\epsilon_{2}}{2} v^{2}, & \sigma > 0, \;\alpha \in \R,\; \epsilon_{i} \in \mathbb{C}\, (i = 1, 2),\\ 
(v(0), u(0)) = (v_0, u_0), 
\end{cases}
\end{equation*}
 posed on a $d$-dimensional  compact Zoll manifold $M$. Considering $\sigma=\frac{\theta}{\beta}$ with $\theta, \beta\in \{n^2:n\in\mathbb{Z}\}$ we derive a bilinear Strichartz type estimate and use it to prove the local well-posedness results for given data $(v_0, u_0)\in H^s(M)\times H^s(M)$ whenever $s>\frac{1}{4}$ when $d = 2$ and $s > \frac{d - 2}{2}$  when $d \geq 3$. Moreover,  in dimensions $2$ and $3$, we use a Gagliardo-Nirenberg type inequality and conservation laws to prove that the local solution can be extended globally in time whenever  $s \geq 1$.
\end{abstract}

\maketitle

\section{Introduction}\label{SectionAC4}
In this work we are interested in addressing  some well-posedness issues to the following  initial value problem (IVP) associated to a system involving  nonlinear Schr\"odinger (NLS) equations with quadratic nonlinearities
\begin{equation}\label{SHGS}
\begin{cases}
i \partial_{t} v \pm \Delta_{g} v -  v = \epsilon_{1} u \bar{v}, & t \in \R,\; x \in M, \\ 
i \sigma \partial_{t} u \pm  \Delta_{g} u - \alpha u = \frac{\epsilon_{2}}{2} v^{2}, & \sigma > 0,\; \alpha \in \R, \;\epsilon_{i} \in \mathbb{C}\; (i = 1, 2),\\ 
(v(0), u(0)) = (v_0, u_0), 
\end{cases}
\end{equation}
 where $v = v(t, x)$ and $u = u(t, x)$ are complex functions, $M$ is a compact Zoll manifold of dimension
 $d \geq 2$ and $\Delta_{g}$ is the Laplace-Beltrami operator.
 
The class of Zoll manifolds is defined as  being that  formed by all compact manifolds such that all geodesics are closed and possess the same period. In particular, this class contains all spheres as well as compact symmetric spaces of rank one (see \cite{BES1978B}, Chapter 4 for a rigorous exposition). The principal reason that  motivated us to consider the system \eqref{SHGS} posed on  Zoll manifolds of dimension
 $d \geq 2$ is that  the spectrum of the Laplacian consists of clusters of
bounded width centered at the points $(k + \frac{Z_0}{4})^2$, $k \in \mathbb{N}$, where $Z_0$ is an integer which depends on the geometry of $M$  (see \cite{WEI1997,ZEL1996,ZEL1997}). Note that the spectrum of the Laplacian on the sphere $\mathbb{S}^d$ 
is exactly of this form (see \eqref{SpecSd} below). This special feature of the spectrum allows  to deduce some arithmetical properties which help us  to understand the behavior of the       Schr\"odinger type groups associated to \eqref{SHGS} for some values of the parameter $\sigma >0$ and consequently to obtain  bilinear estimates that are the main tools to get the well-posedness results (see Propositions \ref{ZollSpec} and \ref{BiliEvo2} below).

The  system of equations \eqref{SHGS} appears in the study of non-linear optics, more specifically
  in studies related to the second harmonic generation (SHG) of type $I$,
 also known as \textit{Frequency Doubling}, which is a nonlinear optical process discovered in the early 1960s. At that time, thanks to the invention of lasers, physicists have obtained a powerful source of coherent light, so that many of the
non-linear optical effects, such as SHG, were demonstrated (see \cite{ZHAOSHI}).
 The functions $v$ and $u$ represent,
respectively, the amplitudes of the envelopes of the first and second harmonics of
an optical wave. For mathematical derivation of similar quadratic models posed on the whole space $\R^n$ and a detailed study of the associated Cauchy problems, we refer to the recent work  \cite{CMS2016}. The SHG system \eqref{SHGS} posed on Riemannian manifold $M$ with the non-Euclidean metric $g$ describes the interaction of these harmonics in a medium in which the optical index is variable. In the system \eqref{SHGS}, there are
four combinations of signs $(\pm, \pm)$ which  are determined by
 the signs of dispersions/diffractions (temporal/spatial cases respectively). The constant
$\sigma > 0 $ measures the dispersion/diffraction rates and plays an important role in the local and global theory. We observe that in order to establish the local well-posedness theory for \eqref{SHGS} it is necessary that $\sigma$  assumes fractional values which is a contrast with the necessary assumption in the one-dimensional case.  The parameter $\alpha > 0$  is 
dimensionless, and in  the one-dimensional physical model one needs to have $0 < \alpha < 1 $ (see \cite{YCM1999}).

If one considers  $\epsilon_{j} \in \{ 1, -1 \}$,  $j = 1, 2$,  for sufficiently regular solutions, the following quantities 
\begin{equation}\label{MASS1}
\mathcal{M}(t) = \int_{M} |v(t)|^{2} + 2 \sigma |u(t)|^{2} dg = \mathcal{M}(0),
\end{equation}
and
\begin{equation}\label{ENERGY2}
\mathcal{E}(t) = \int_{M} \left( |\nabla v(t)|_{g}^{2} +  |\nabla u(t)|_{g}^{2} + |v(t)|^{2} + \alpha |u(t)|^{2} + \epsilon_{1} \Ree(v^{2}(t) \overline{u(t)})\right) dg = \mathcal{E}(0),
\end{equation}
where $\Ree$ denotes the real part, are conserved. The quantities in \eqref{MASS1} and \eqref{ENERGY2} represent respectively the mass and energy of the system \eqref{SHGS}.

For simplicity of exposition we consider  the combination  $(+, +)$ of signs and use the notation $\Delta$ to represent $\Delta_g$ throughout this work. The study of the Schr\"odinger equations 
with quadratic non-linearities has attracted attention of several mathematicians over the past decades, see for instance \cite{TB2006,KPV96,LH2014} and references therein. As far as we know, in the literature, the system \eqref{SHGS} has been studied considering
$M = \mathbb{R}^{d}$,  $(1 \leq d \leq 5)$ in \cite{ZHAOSHI} where the focus is on the variational questions, and  $ M = \mathbb {T} $ in  \cite {AL2007} where the focus is on the  local well-posedness issues for $(u_0, v_0) \in H^{s}(\mathbb{T})\times H^{s}(\mathbb{T})$. More precisely, in \cite {AL2007}  the local well-posedness results are obtained for given data with regularity $ s \geq 0 $ if $1 / \sigma > 0$
and $ s > -1/2 $ if $ 1 / \sigma \neq 1 $. In addition, the authors in \cite{AL2007} used  mass conservation \eqref{MASS1} and proved  global well-posedness  for $ s \geq 0 $. Also, 
we can cite \cite{HLO2011} where the questions about scattering theory in $ \R^d, d  \geq 3 $ are addressed for a similar system  to \eqref {SHGS}, 
and \cite {HOT2013} where a study of well-posedness and blow-up is performed. 

As mentioned earlier, if we choose  $(+, +)$  signs and $\sigma> 0 $,  the system \eqref{SHGS} can be rewritten 
as
\begin{equation}\label{SHGSA}
\begin{cases}
i \partial_{t} v + \Delta v -  v = \epsilon_{1} u \bar{v} \\ 
i  \partial_{t} u +  \frac{1}{\sigma} \Delta u - \frac{\alpha}{\sigma}  u = \frac{\epsilon_{2}}{2 \sigma} v^{2} \\ 
(v(0), u(0)) = (v_{0}, u_{0}).
\end{cases}
\end{equation}
 The main  objective of this work is in addressing the well-posedness issues for the system \eqref{SHGSA} posed on $d$-dimensional compact Zoll manifold $M$ with  given  data $(u_0, v_0)$
in a suitable Sobolev spaces $H^s(M)\times H^s(M)$. As far as we know,  the results involving   the well-posedness theory  for  \eqref {SHGSA} posed on compact manifolds in dimension $d \geq 2$ are not known.
 
To accomplish our objective, using Duhamel's formula, we can consider   the IVP \eqref{SHGSA} in the  following equivalent integral  formulation 
\begin{equation}\label{SHGSB}
\begin{cases}
v(t):= V(t)v_0 - i \epsilon_1  \displaystyle \int_{0}^{t} V(t - t') \overline{v}(t') u(t') dt',\\
u(t):= U_{\sigma}(t)u_0 - i  \frac{\epsilon_2}{2 \sigma} \displaystyle \int_{0}^{t}  U_{\sigma}(t - t')v^{2}(t') dt',
\end{cases}
\end{equation}
where $V(t) = e^{i t (\Delta - 1)}$ and  $U_{\sigma}(t) =  e^{i t (\frac{1}{\sigma}\Delta - \frac{\alpha}{\sigma})} $ 
are the respective unitary groups associated with the linear problem. From now on, we consider the equivalent system \eqref{SHGSB} and use the contraction mapping argument in an appropriate space to get the required solution. Before announcing  the main  results on well-posedness theory, we introduce some definitions and function spaces on which we will be working.  For convenience, let us start recapitulating some notions on Sobolev spaces in compact manifolds.

 We denote by
  $\{e_{k}\}_{k \geq 0} \subset L^{2}(M, \mathbb{C})$ 
an orthonormal basis formed by eigenfunctions of   $-\Delta$, with eigenvalues $ \{\mu_{k}\}_{k \geq 0}$ and by $P_{k} : L^2(M) \rightarrow L^2(M)$ the orthogonal projection
on $e_{k}$, given by $$ P_{k} f =  \langle e_{k}, f \rangle_{L^{2}(M)} e_{k}.$$ In this way, we can define the Sobolev space $H^{s}(M)$ as being the completion of the space $C_0^{\infty}(M)$ with respect to the norm 
\begin{equation}\label{D2}
\|u\|_{H^{s}(M)}^{2} :=   \sum_{k \geq 0} \langle \mu_{k}\rangle^{s} \|P_{k} u\|_{L^{2}(M)}^{2} \simeq \|(1 - \Delta)^{\frac{s}{2}} u\|_{L^{2}(M)}^{2}.
\end{equation}

%%%%%%%%%%%%%%%%%%%%%%%%%%%%%%%%%%%%%%%%%%%%%%%%%%%%%%%%%%%%%%%%%%%%%%%%%%
In what follows, we introduce a generalization of the spaces $X^{s,b}(\R \times \mathbb{T}^{d})$ previously  introduced by Bourgain in \cite{BOU93A} in the context of the NLS equation. The generalization in the context  of the compact manifolds $M$  is due to
Burq, G\'erard and Tzvetkov in \cite{BGT1,BGT2,BGT3} where the authors studied  the NLS equation. The definition  uses the structure of the spectrum of  $- \Delta$. 

 \begin{defini} \label{D2a}
Let $s, b \in \mathbb{R}$. The space $X^{s,b}(\mathbb{R}\times M)$ is the completion of the space  $C_{0}^{\infty}(\mathbb{R} \times M)$ with respect to the norm
\begin{equation}\label{D3}
 \|u\|^{2}_{X_{- \Delta}^{s,b}(\mathbb{R} \times M)} = \sum_{k} \| \langle \tau  + \mu_{k} \rangle^{b} \langle \mu_{k}\rangle^{\frac{s}{2}} \widehat{P_{k} u}(\tau)\|_{L^{2}(\mathbb{R};L^{2}(M))}^{2}
= \| S(-t)  u(t, \cdot)\|_{H^{b}(\mathbb R_{t}; H_{x}^{s}(M))}^{2},
\end{equation}
where $S(-t) := e^{- i t \Delta}$  and  $\widehat{P_{k} u}(\tau)$ denotes the Fourier transform of the function $t \mapsto P_{k} u(t, \cdot)$. 
\end{defini}

Taking  in consideration the  spaces $X^{s,b}$ given by  Definition \ref{D2a}, we introduce a family of spaces associated with the linear  structure of \eqref{SHGSA}. Our definition is appropriate to the
 modulation produced by the unitary groups $V (t)$ and $U_{\sigma}(t)$.

\begin{defini}\label{D2aa} Let $ \delta, \gamma  \in \mathbb{R}$ be fixed. Given $s, b \in \mathbb{R}$, 
we define the spaces $X^{s,b}_{\delta, \gamma}(\mathbb{R}\times M)$, as being the completion of the space  $C_{0}^{\infty}(\mathbb{R} \times M)$ with respect to the norm
\begin{equation}\label{D3bb}
\begin{split}
 \|u\|^{2}_{X_{- \delta  \Delta + \gamma }^{s,b}(\mathbb{R} \times M)} &= \sum_{k} \Big\| \langle \tau  + \delta  \mu_{k} + \gamma \rangle^{b} \langle \mu_{k}\rangle^{\frac{s}{2}} \widehat{P_{k} u}(\tau) \Big\|_{L^{2}(\mathbb{R};L^{2}(M))}^{2} \\
 & = \sum_{k}  \langle \mu_{k}\rangle^{s} \int_{\R} \langle \tau  + \delta  \mu_{k} + \gamma \rangle^{2b} \| \widehat{P_{k} u}(\tau)\|_{L^{2}(M)}^{2} d \tau \\
&= \left\| e^{ it (-\delta \Delta + \gamma )  } u(t, \cdot) \right\|_{H^{b}(\mathbb R_{t}; H_{x}^{s}(M))}^{2}.
\end{split}
\end{equation}
 In particular, 
if $\delta \neq 0 $ and $\gamma =0$,  we denote this space by  $X^{s,b}_{\delta}$.
\end{defini}

In order to establish the local theory, we need to define a local version of the spaces
 $ X^{s,b}_{ \delta,  \gamma} (\mathbb{R} \times M)$ with respect to the variable $t \in I \subset \R $. 
 
\begin{defini}\label{D5}
Let  $I \subset \mathbb{R}$ be a compact interval. We define the restriction space $X_{\delta, \gamma}^{s,b}(I \times M)$  equipped with the following norm
\[
\|u\|_{X_{\delta, \gamma}^{s,b}(I \times M)} := \inf_{w \in X_{\delta, \gamma}^{s,b}(\mathbb{R} \times M)} \Big\{ \|w\|_{X_{\delta, \gamma}^{s,b}(\mathbb{R} \times M)} : w\big|_{I}\, = u \Big\}.
\]
\end{defini}

\begin{obs}\label{XsbObs}
The definition (and the norm) of the spaces  $X^{s,b}_{\delta, \gamma}$ clearly depend on the operator 
$P := - \delta  \Delta + \gamma $. However, if there is an operator $ Q $ of the same order as $ P $, having the same eigenfunctions and such that the eigenvalues $p_{k}$ and  $q_{k}$ of  $P$ and $Q$ 
respectively, obey the condition
$$|p_{k} - q_{k}| \leq C_{0}, $$ 
for some constant $C_{0} > 0$, then one can easily  show that  
there exists $C > 0$ such that for all $k \in \mathbb{N}$ and  $\tau \in \mathbb{R}$, 
\begin{equation*}\label{x1}
\frac{1}{C} \langle \tau + p_{k} \rangle \leq \langle \tau + q_{k} \rangle 
\leq C \langle \tau + q_{k} \rangle .
\end{equation*}
Consequently, $X^{s,b}_{P}$ and $X^{s,b}_{Q}$ have equivalent norms.

In particular, the spaces 
$X^{s,b}_{- \delta  \Delta + \gamma} $  and $ X^{s,b}_{- \delta  \Delta }$, have equivalent norms, 
\[
\| \cdot \|_{X^{s,b}_{- \delta  \Delta + \gamma}}  \simeq   \| \cdot \|_{X^{s,b}_{- \delta  \Delta}} =:\| \cdot \|_{X^{s,b}_{\delta}}.
\]
\end{obs}

Hence, according to the norm equivalence given by Remark  \ref{XsbObs},  instead of $X^{s,b}_{\delta, \gamma}$
we can work in the $X^{s,b}_{\delta }$ and  $X^{s,b}$ spaces.
%%%%%%%%%%%%%%%%%%%%%%%%%%%%%%%%%%%%%%%%%%%%%%%%%%%%%%%%%%%%%%%%%%%%%%%%%%
Note that, considering the system \eqref{SHGSB}, to obtain the local well-posedness results using the contraction mapping argument, one needs to establish the following  crucial bilinear estimates
\begin{equation*}\label{BiIntroA}
\|u_{1} \overline{u_{2}}\|_{X_{- \Delta + 1}^{s, -b'}} \leq C \|u_{1}\|_{X_{- \frac{1}{\sigma} \Delta + \frac{\alpha}{\sigma}}^{s, b}} \|u_{2}\|_{X_{- \Delta + 1}^{s, b}}, 
\end{equation*}
and
\begin{equation*}\label{BiIntroB}
\|v_{1} v_{2}\|_{X_{- \frac{1}{\sigma} \Delta + \frac{\alpha}{\sigma}}^{s, -b'}} \leq C \|v_{1}\|_{X_{- \Delta + 1}^{s, b}} \|v_{2}\|_{X_{- \Delta + 1}^{s, b}},
\end{equation*}
for some $C > 0$ and $(b, b') \in \R^{2}$ satisfying $b + b' < 1$ and $0 < b' < 1/2 < b$.  

In view of  Remark \ref{XsbObs}, to obtain the above estimates, it is sufficient to prove that 
\begin{equation}\label{BiIntroAA}
\|u_{1} \overline{u_{2}}\|_{X^{s, -b'}} \leq C \|u_{1}\|_{X_{1/\sigma}^{s, b}} \|u_{2}\|_{X^{s, b}}, 
\end{equation}
and 
\begin{equation}\label{BiIntroBB}
\|v_{1} v_{2}\|_{X_{1/\sigma }^{s, -b'}} \leq C \|v_{1}\|_{X^{s, b}} \|v_{2}\|_{X^{s, b}}
\end{equation}
hold for some $C > 0$ and $(b, b') \in \R^{2}$ satisfying $b + b' < 1$ and $0 < b' < 1/2 < b$.

To prove  \eqref{BiIntroAA} and  \eqref{BiIntroBB} we use  a duality argument  followed by   dyadic decompositions on the functions  $u_j$ and $v_j$, $j=1,2$. The crucial fact in this process is to analyse the decompositions of $ u_1 $ and $u_2$ when they are localized  on the dyadic  frequencies $\sim N$ and $\sim L$  respectively, with $N \ll L$ or $L \ll N$. 

 In view of the equivalence that will be given below in Lemma \ref{Eq00}, we find that such estimates
are closely related with the bilinear Strichartz estimates where discrepancy between frequencies is ``controlled" \;
by  $(\min(N, L))^{s}$, where $s > s_{0}(M)$. More precisely, we will show that
\begin{equation}\label{BilStrA}
\|\tilde{u_{1}} u_{2}\|_{L^{2}(\mathbb{R} \times M)} \leq C (\min(N, L))^{s} \|u_{1}\|_{X^{0, b}_{1/\sigma}(\mathbb{R} \times M)} 
 \|u_{2}\|_{X^{0, b}(\mathbb{R} \times M)}, 
\end{equation}
with $\tilde{u_1} = u_{1}$ or $\overline{u_1}$, is equivalent to 
\begin{equation}\label{BilStrB}
    \|e^{ \pm i  \frac{t}{\sigma}\Delta}u_0 \; e^{i t \Delta} v_0\|_{L^{2}((0,1)_{t} \times M)} \leq C (\min(N, L))^{s} \|u_0\|_{L^{2}(M)}\|v_0\|_{L^{2}(M)},
\end{equation}
where $ N, L $ are dyadic numbers on which  $ u_0 $ and $ v_0 $ are spectrally localized respectively. In the next section we highlight some difficulties in
the proof of the bilinear estimate \eqref{BilStrB}. 

\section{Arithmetical properties of the spectrum and applications to bilinear Strichartz estimates}\label{Section2}

 One of the useful properties of the Laplace-Beltrami operator on compact manifolds  is that 
it has a discrete spectrum. In the case $M = \mathbb{T}^2$, this property  as well as  expansion in Fourier series  is widely used. For instance, in Bourgain's pioneer works \cite{BOU93A,BOU93B} the $L^4$--Strichartz estimate was reduced to counting the number of elements in the set 
\[
A(k) = \{ (n_1,n_2) \in \mathbb{Z}^2 ; |n_1|, |n_2| \leq N, n_1^2 + n_2^2 = k \}.
\]
It is well known that,  $\sharp A(k) \lesssim N^{\epsilon} $  for any $\epsilon >0$. 
However, in \eqref{BilStrB} the product of Schr\"odinger semigroups with distinct
spectra  produces a different dispersion phenomena, especially due to the constant $\sigma > 0$. Bilinear estimate to spectral projectors was proved by  Burq, G\'erard, Tzvetkov in \cite{BGT2} for the case of compact surfaces (see Proposition \ref{BiliSurf} below) and in \cite{BGT1, BGT3} for  the  higher dimensional case, viz., $\mathbb{S}^{d}$ $d \geq 2$  (see Proposition \ref{BiliSd} below). 
 In general, the knowledge of the spectrum is combined with some estimates from the analytic number theory. At this point, we highlight that, as  in \cite{BOU93A, BGT1, BGT2}, estimates on the upper bound for the number of intersections of a closed curve  with a lattice are strongly used. For example, an estimate  for the number of intersection points of oval shaped closed curves with  $ \mathbb{Z}^ 2$, is following result  due to Bombieri, Pila \cite {BP1989}.

\begin{teorema}\label{BPTeo}
Suppose $\phi: \mathbb{S}^{1} \rightarrow \mathbb{R}^{2}$ analytic. Then, for all $\varepsilon >0$,
\[
|t \phi(\mathbb{S}^{1}) \cap \mathbb{Z}^{2}| \leq C(\phi, \varepsilon) t^{\varepsilon}.
\]
\end{teorema}
Besides, it is proved that if $C$ is an affine plane algebraic curve of degree $d$ with integer
coefficients, and  if one takes a box with
sides of length $N$, then $C$ can contain no more than $ O_{d, \epsilon}(N^{\frac{1}{d}+ \varepsilon})$
integer points
within the box. 

Nevertheless, in the case of the plane algebraic curves of the form   $ \{(x, y) \in \mathbb{R}^2: a x^{2} - b y^{2} = m \} $ where $m \in \mathbb{Z}$ and $a, b \in \R$, Theorem \ref{BPTeo} is not valid. 

% In our case too, we need to introduce  a similar constraint (see Proposition  \ref{BiliEvo}).

 To describe better the situation which appears in our analysis, let us consider $ M =\mathbb{S}^2$. Suppose that $N, L$ are natural numbers with  $N \leq L$. Due to the fact that $ Spec(- \Delta_{\mathbb{S}^2}) = \{ j (j + 1) \}_{j \in \mathbb{N}}$,  we need to prove that   the number of solutions of the inequality 
\begin{equation}\label{SigmEq}
\left|k^2 + k -  \sigma (\ell^2 + \ell) - \sigma m \right| < \frac{\sigma}{2}
\end{equation}
where $N \leq k \leq  2 N$ and $m \in \mathbb{Z}$ is  $O(N^{\varepsilon})$. In general, we don't know how to estimate the number of solutions to  \eqref{SigmEq}. However, if we  choose the parameter $\sigma>0$ in the following form
\begin{equation}\label{restriction} 
\sigma = \frac{\beta}{\theta}\qquad \textrm{ with}\quad \beta, \theta \in \{ n^{2}: n \in \N \},
\end{equation}
 we obtain that \eqref{SigmEq} is reduced to the problem of find the number of solutions of the equation 
\begin{equation}\label{SigmEq-1}
 (\sqrt{\theta} (2 k + 1))^{2} - (\sqrt{\beta} (2 \ell + 1))^{2} = p,
\end{equation}
for some explicit values of $p \in \mathbb{Z}$. Thus, to consider this case,  we will use Lemma \ref{NTLema} to estimate the number of solutions to \eqref{SigmEq-1}. This is the point, that lead us to impose restriction \eqref{restriction} on the parameter $\sigma$.
It is worth mentioning that:

\begin{enumerate}
    \item[$\bullet$]  Using the Strichartz estimates proved by Burq, G\'erard and Tzvetkov in \cite{BGT3}, we may obtain  local well posedness for the quadratic system \eqref{SHGS} without any restriction on the parameter $\sigma >0$ for general compact manifolds
    of dimension $d \geq 2$ for initial data in $\bH^{s}(M)$ with $s > \frac{d - 1}{2}$, by using a procedure as in  Nogueira and Panthee \cite{NP2020}.  However, in our case, i.e., in the case of compact Zoll manifolds, the use of the  Bourgain's spaces improves the range of Sobolev index for $s > \frac{1}{4}$ for $d = 2$  and $s > \frac{d -2}{2}$ for $d \geq 3$. In particular, in dimension 3, this 
    allows us to extend the solutions (with initial data in $\bH^{1}(M^3)$ ) globally in time, 
    which is not possible in the case of general $3$-dimensional compact manifolds in the case of the cubic nonlinear Schr\"odinger equation (see  \cite{BGT3} p. 571).
    
    \item[$\bullet$] A result due to Huxley \cite{HUX03}
    asserts that if $Q(x,y) = a x^2 + b xy + c y^2$ is a quadratic form with 
    $a > 0$ and $D:= 4 a c - b^2 > 0$, then 
    \[
    \sharp \{ (x,y)  \in  \mathbb{Z}^2 : Q(x,y) \leq N\} = \frac{2 \pi}{\sqrt{D}} N + O(N^{\frac{131}{416}+}),
    \]
    when $N \rightarrow + \infty$. Observe that considering \eqref{SigmEq} and completing squares, we can suppose that the associated quadratic form is $Q(X,Y) = \frac{1}{4} X^2 -  \frac{\sigma}{4} Y^2$, where we performed the change of variables $X = 2 x + 1$ and $Y = 2y + 1$. Thus, since $b = 0$, we have $D = - \frac{\sigma}{4} < 0$. So, the above estimate is not applicable in our case.  However, this result is very important in the study of the local well-posedness for  the NLS equation on the bidimensional irrational tori, see Demirbas \cite{DEM2017}. 
    
     \item[$\bullet$] The recent methods via the application of the Bourgain-Demeter \cite{BD1,DEM18} decoupling theory seem promising to solve many questions related to questions of this nature. For example, let $u_{0} \in \mathbb{T}^{d -1}$ (flat torus), with $\supp(\widehat{u_0}) \subset [-N, N]$. Then, for each $\varepsilon >0$ an application of the decoupling theory in \cite{BD1} guarantees  the full range of expected $L_{x t}^p$ Strichartz estimates
     \[
     \|e^{i t \Delta_{\mathbb{T}^{d - 1}} } u_0\|_{L^{p}([0, 1] \times \mathbb{T}^{d - 1})} \lesssim_{\varepsilon} N^{\varepsilon} \|u_0\|_{L^2}
     \]
     if $2 \leq p \leq \frac{2 (d + 1)}{d - 1}$, and 
      \[
     \|e^{i t \Delta_{\mathbb{T}^{d - 1}} } u_0\|_{L^{p}([0, 1] \times \mathbb{T}^{d - 1})} \lesssim_{\varepsilon} N^{\frac{d - 1}{2} - \frac{d+ 1}{p} + \varepsilon} \|u_0\|_{L^2}
     \]
     if $p \geq \frac{2 (d + 1)}{d - 1}$. For more details we refer to the recent work of Demeter \cite{DEM20}.
     
     Recently,  Fan et. al \cite{FSWW2018} used the decoupling type argument from \cite{BD1} to obtain the bilinear Strichartz type estimate for irrational tori,  recovering and generalizing the result of  \cite{SPST07}.

     The decoupling theory methods are applied to obtain estimates for  exponential sums whose phase function is associated with some curved surface, for example a truncated cone. 
     However,  we do not have idea how to use these methods to improve the estimates for our bilinear Strichartz estimates, which are of the form 
     \[
    \|e^{ \pm i \frac{t}{\sigma} \Delta_{M}} u_0 e^{i t \Delta_{M}}v_0 \|_{L^{p}([0,1] \times M)}  
    \]
   % \[
    %= \Big\| \sum_{k} \sum_{\ell} a_k a_{\ell} e^{i t (\frac{k^2 + k}{\sigma})} e^{- i t (\ell^2 + \ell)} \Big\|_{L^{p}([0,1] \times M)}.
    % \]
     where $M$ is a Zoll manifold or the $d$-dimension sphere and $u_0, v_0 \in L^2(M)$, 
     in the general case $ \sigma > 0 $, i.e., without any restriction on $\sigma$.
\end{enumerate}

\begin{obs}
  This sort of problems are also studied  for algebraic plane curves of varying  degrees (see \cite{EV05,VAU14}). 
\end{obs}
%%%%%%%%%%%%%%%%%%%%%%%%%%%%%%%%%%%%%%%%%%%%%%%%%%%%%%%%%%%%%%%

\section{Main Results}
 In this section,  we  state the main results on the  well-posedness for the IVP \eqref{SHGSA} posed on the class of $d$-dimensional Zoll manifolds $M$ that we  introduced in the beginning of the introduction. 
To simplify the notations, we define $\bH^{s}(M):=H^{s}(M)\times H^{s}(M)$ and consider the IVP \eqref{SHGSA} with initial data $(u_0, v_0) \in \bH^{s}(M)$.

Before stating the main  results, let $d \geq 2$ and  define 
\begin{equation}\label{BiliSdb}
s_{0}(d) :=\begin{cases}
\frac{1}{4}   &\mbox{ if } d = 2  \\
\frac{1}{2} +  \varepsilon   &\mbox{ if } d = 3 \\
 \frac{d-2}{2}    &\mbox{ if } d \geq 4.
  \end{cases}
\end{equation}
The  main result concerning the local well-posedness theory  in $\bH^{s}(M)$ for $s > s_{0}(d)$ is the following.

\begin{teorema}\label{BCL00}
Let  $(M, g)$ be a $d$-dimensional Zoll manifold  and $\sigma = \frac{\beta}{\theta}$ with $\beta, \theta \in \{ n^{2}: n \in \N \}$. For any  $(v_0,u_0) \in \bH^{s}(M)$, with 
  $ s > s_0(d)$ where $s_0(d)$ is defined in \eqref{BiliSdb}, there exist  $T = T(\|(v_0,u_0)\|_{\bH^{s}}) > 0$
and a unique solution  $(v(t), u(t))$ of the IVP \eqref{SHGSA} on the interval $[0,T]$ such that, for some   $b > \frac{1}{2}$
\begin{enumerate}
\item[$(i)$] $(v, u) \in X_{T}^{s,b}(M) \times (X_{1/\sigma}^{s, b})_{T}(M)$,
\item[$(ii)$] $(v,u) \in C([0, T], \bH^{s}(M))$.
\end{enumerate}
 
 Moreover, for any $T' \in (0, T)$ there exists an open ball $B((v_0, u_0), r) \subset \bH^{s}(M)$ such that the application   $$ \Phi: B((v_0,u_0), r) \ni (\tilde{v_0}, \tilde{u_0}) \mapsto (\tilde{v}, \tilde{u}) \in C([0, T'], \bH^{s}(M)) ,$$ is   Lipschitz-continuous
 for some $r >0$. 
\end{teorema}

%%%%%%%%%%%%%%%%%%%%%%%%%%%%
\begin{obs}
A natural question is whether the local result given by Theorem \ref{BCL00} is sharp. As far as we know, even for $d=2$, the local well-posedness for the cubic or quadratic nonlinearities for  $s = 1/4$ in classical  Sobolev spaces $H^s(M)$  are open problems in the context of compact manifolds. However, using Besov spaces, Takaoka \cite{TAK2016} proved
that the  local well-posedness can be  achieved for $s = 1/4$ for $M = \mathbb{S}^2$. The analogous question  for the quadratic system \eqref{SHGS}
is being addressed in a work   in progress by the first author in  \cite{MN202X}. Ill-posedness results will be addressed elsewhere.
\end{obs}
%%%%%%%%%%%%%%%%%%%%%%%%%%%

 Note that the system \eqref{SHGSA} has mass and energy conservation laws given respectively by \eqref{MASS1} and
  \eqref{ENERGY2}. For $s \geq 1$ we use these conserved quantities together with a Gagliardo-Nirenberg inequality (see Proposition \ref{GNM} below) to prove that the solutions given by Theorem \ref{BCL00}
  in the case $ d = 2, 3$  are in fact global. This is the content of the following theorem.

\begin{teorema}\label{BCG00}
Let $(M, g)$ be a $d$-dimensional Zoll manifold, $s \geq 1$ and  $T^{\ast}  > 0$  be the maximal time of existence for the solution
\[
(v, u) \in C([0, T^{\ast}), \bH^{s}(M))
\]
given by  Theorem \ref{BCL00} when $d = 2, 3$. Then
$T^{\ast} = + \infty$, that is, the solution is global in time. 
\end{teorema}

This paper is organized as follows. In Section \ref{SectionBC4}  we record some basic properties of the  $ X^{s,b}$ spaces. In  Sections \ref{SectionCC4} and  \ref{SectionDC4} we will derive the  bilinear estimates that are crucial in our argument. The proofs of the main results are presented in Section \ref{SectionEC4}  with subsection \ref{SSE2C4} devoted to supply the proof for Theorem \ref{BCL00}, and subsection  \ref{SSE3C4}  for Theorem \ref{BCG00}. 

%%%%%%%%%%%%%%%%%%%%%%%%%%%%%%%%%%%%%%%%%%%%%%%%%%%%%%%%%%%%%%%%
\section{ Basic properties of the  Bourgain spaces} \label{SectionBC4}
%%%%%%%%%%%%%%%%%%%%%%%%%%%%%%%%%%%%%%%%%%%%%%%%%%%%%%%%%%%%%%%%%
We begin exploring some basics properties of the Bourgain  spaces $X^{s,b}$ and  $X^{s,b}_{\delta}$. Here, $M$ denotes a general $d$-dimensional compact manifold. 

%In view of the norm equivalence given by Remark \ref{XsbObs}, 
%we can prove the statements in this section with $\delta = 1$, 
%and consider the spaces $X^{s, b}$ rather than $X^{s,b}_{\delta}$. 

\begin{prop}\label{BasicXsb} The following properties are valid 
\begin{enumerate}
\item[(i)] For  $s_{1} \leq s_{2}$ and $b_{1} \leq b_{2}$, one has  $X^{s_{2},b_{2}}(\mathbb{R}\times M) \hookrightarrow X^{s_{1},b_{1}}(\mathbb{R}\times M)$.
\item[(ii)] $X^{0,\frac{1}{6}}(\mathbb{R}\times M) \hookrightarrow L^{3}(\mathbb{R}, L^{2}(M))$. 
\item[(iii)] If  $b > \frac{1}{2}$,  then the inclusion $X^{s,b}(\mathbb{R} \times M) \hookrightarrow C(\mathbb{R}, H^{s}(M))$ holds.
\end{enumerate}
\end{prop}
\noindent
\begin{dem}The part $(i)$ follows directly from the Definition \ref{D2a}. The part $(ii)$ follows
from the fact that $u \in X^{s, b}(\mathbb{R} \times M)$, if  and only if,  $S(-t) u (t, \cdot) \in H^{b}(\mathbb{R}, H^{s}(M))$
and from the immersion $H^{1/6}(\mathbb{R}) \hookrightarrow L^{3}(\mathbb{R}) $,  (see \cite{ANT2008A} p. 1866, for more details).
The proof of $(iii)$, is given in \cite{BGT2}, p. 196.
\end{dem}

\begin{prop}\label{Eq01}
Let  $ u \in C_{0}^{\infty}(\mathbb{R} \times M)$ and  $(s, b) \in \R^{2}$.  Considering the expression 
\begin{equation}\label{Eq1}
w(t) = \frac{1}{2\pi} \sum _{k \in \mathbb{N}} \langle \mu_{k}\rangle^{\frac{s}{2}} \int_{\mathbb{R}}  \langle \tau + \mu_{k}\rangle^{b} \widehat{P_{k} u}(\tau) e^{i t \tau} d\tau,
\end{equation}
one has 
\begin{equation}\label{Eq2}
u(t) = \frac{1}{2\pi} \sum _{k \in \mathbb{N}} \langle \mu_{k}\rangle^{-\frac{s}{2}} \int_{\mathbb{R}}  \langle \tau + \mu_{k}\rangle^{-b} \widehat{P_{k} w}(\tau) e^{i t \tau} d\tau.
\end{equation}
Moreover, 
\begin{equation}\label{NormEquality}
\|w \|_{L^{2}(\mathbb{R} \times M)}  =  \|u\|_{X^{s,b}(\mathbb{R} \times M)}.
\end{equation}
\end{prop}
\noindent
\begin{prova}
 The proof follows by applying  orthogonal projection, Fubini's theorem and  Planchrel's identity.
\end{prova}

We note that the results stated in Propositions \ref{BasicXsb} and \ref{Eq01} hold for the $X^{s,b}_{\delta}$ spaces as well.

%%%%%%%%%%%%%%%%%%%%%%%%%%%%%%%%%%%%%%%%%%%%%%%%%%%%
\subsection{Spectral projectors}\label{projetors}
%%%%%%%%%%%%%%%%%%%%%%%%%%%%%%%%%%%%%%%%%%%%%%%%%%%%%

In this section, we will introduce spectral projection operators and their properties which are  used in the decomposition of functions. Also we record some estimates involving these operators that are used  especially in the following section, where we derive some estimates involving the spectral localization of the functions related to the  $X_{\delta, \gamma}^{s,b}$ spaces.
In what follows, we use $N = 2^{n}, L = 2 ^{\ell}$, $n, \ell \in \mathbb{N}$, to denote the dyadic integers. For clarity of exposition we consider $\delta = 1$ and $\gamma = 0$, the similar results hold for $\delta \ne1$ and $\gamma \ne 0$  too. 

\begin{defini}
Given a dyadic integer $N$, we say that a function $u$ on $M$ 
is spectrally localized at frequency $N$ if 
\begin{equation}\label{LocU}
u = u^{N} = \Pi_{N}(u) := \sum_{ k: N \leq \langle \mu_{k} \rangle^{1/2} < 2N} P_k u. 
\end{equation}

Moreover, we denote $\Pi_1 =  \sum_{ k: 0 \leq \langle \mu_{k} \rangle^{1/2} < 2} P_k u$.
\end{defini}

Let  $u = u(t, x) \in C_{0}^{\infty}(\mathbb{R} \times M)$.  Using a  $L^{2}$ -decomposition  of $u(t, \cdot)$, we can write 
\[
u = \sum_{k} P_{k}u =\sum_{N \geq 1} \Pi_{N}(u).
\]
Now, applying the Fourier inversion theorem 
to the function $t \mapsto \Pi_{N} u(t, \cdot)$, we obtain
\[
u  =\sum_{N \geq 1} \left( \frac{1}{2 \pi} \int_{\mathbb{R}}\widehat{\Pi_{N} u} (\tau) e^{it \tau} d\tau \right)  =  \sum_{N \geq 1} \sum_{L \geq 1} \Pi_{NL}(u),
\]
where
\begin{equation}\label{PiNL}
 \Pi_{NL}(u) :=   \sum_{k:  N \leq\langle \mu_{k} \rangle^{1/2} < 2 N}  \frac{1}{2 \pi} \int_{L \leq \langle \tau + \mu_{k} \rangle \leq 2 L}\widehat{P_{k} u} (\tau) e^{it \tau} d\tau =: \sum_{k:  N \leq\langle \mu_{k} \rangle^{1/2} < 2 N} \Pi_{k,L}(u).
\end{equation}

According to this construction, we can identify two types of localization operators, the one with the spatial variable and the next with  the time  variable. More precisely, we have
\begin{enumerate}
\item[$\bullet$] Localization with respect to  the time variable: 
\begin{equation}\label{uL}
 u^{L} = \displaystyle \textbf{1}_{L \leq\langle \tau + \mu_{k} \rangle \leq 2 L} (u) = \Pi_{L}(u):= \frac{1}{2 \pi}\int_{\Lambda _{L}} \widehat{u} (\tau) e^{it \tau} d\tau,
\end{equation}
where  $ \Lambda_{L} := \{ \tau \in \mathbb{R} ; L \leq\langle \tau + \mu_{k}\rangle \leq 2 L \}$.

\item[$\bullet$] Localization with respect to the space variable:
\begin{equation}\label{uN}
u^{N}=\Pi_{N}(u):= \sum_{k:  N \leq\langle \mu_{k} \rangle^{1/2} < 2 N} P_{k}u .
\end{equation}

\end{enumerate}

%From  \eqref{PiNL}, denoting $\Pi_{NL}(u)$ by $u^{NL}$, we can write
%\begin{equation}\label{uNL}
% u^{NL}
% =  \frac{1}{2 \pi}   \int_{L \leq \langle \tau + \mu_{k} \rangle \leq 2 L}  \widehat{ u^{N}} (\tau) e^{it \tau} d\tau.
%\end{equation}
%This shows that   $u^{NL}$ can be viewed as  $(u^{N})^{L}$. Thus
%\begin{equation*}
%\begin{split}
%\textbf{1}_{L \leq\langle \tau + \mu_{k} \rangle \leq 2 L}(\Pi_N(u)) =\frac{1}{2\pi} \int_{\Lambda_{L}} \sum_{k:  N \leq\langle \mu_{k} \rangle^{1/2} < 2 N} \widehat{P_{k} u}(\tau) e^{it \tau} d \tau = u^{NL}.
%\end{split}
%\end{equation*}
%On the other hand, it follows from definitions and the  Fubini's  theorem, that
%\begin{equation*}
%\Pi_{N} (\textbf{1}_{\langle \tau + \mu_{k}\rangle \sim L} u ) =  \frac{1}{2 \pi}  \sum_{k:  N \leq\langle \mu_{k} \rangle^{1/2} < 2 N} \int_{\Lambda_{L}} \widehat{  P_{k}(u)}(\tau) e^{i t \tau } d \tau = u^{NL}.
%\end{equation*}
From the previous definitions, it is easy to see that the operators   $\textbf{1}_{L \leq\langle \tau + \mu_{k} \rangle \leq 2 L}$ and $ \Pi_{N}$  commute with each other and can be applied  in different orders if necessary. 

Using the expression of the projector $\Pi_{NL}$ given in      \eqref{PiNL} and the definition of the $ X^{s,b}$-norm 
 given in Definition \ref{D2a}, we can prove the following basic estimates.

\begin{lema}\label{LemaApB}

Let $s, b \in \mathbb{R}$. Then, there exists  $C > 0$ such that: 
\begin{equation}\label{LemaApB1}
\frac{1}{C} \|u^{NL}\|_{X^{s,b}(\mathbb{R} \times M)} \leq L^{b} N^{s} \|u^{NL}\|_{L^{2}(\mathbb{R} \times M)} \leq C  \|u^{NL}\|_{X^{s,b}(\mathbb{R} \times M)},
\end{equation}
\begin{equation}\label{LemaApB2}
\frac{1}{C} \sum_{N, L} L^{2b} N^{2s} \|u^{NL}\|_{L^{2}(\mathbb{R} \times M)}^{2} \leq \|u\|_{X^{s,b}(\mathbb{R} \times M)}^{2} \leq C \sum_{N, L} L^{2b} N^{2s} \|u^{NL}\|_{L^{2}(\mathbb{R} \times M)}^{2}, 
\end{equation}
where the summation is taken over all dyadic values of 
 $N$ and $L$.  
 
 Similar estimates also hold for $X^{s,b}_{\delta}$.
\end{lema}

%%%%%%%%%%%%%%%%%%%%%%%%%%%%%%%%%%%%%%%%%%%%%%%%%%%%%%%%%%%%%%%%%
\section{ Bilinear estimates and applications}\label{SectionCC4}
%%%%%%%%%%%%%%%%%%%%%%%%%%%%%%%%%%%%%%%%%%%%%%%%%%%%%%%%

In this section, we derive bilinear interaction estimates for the quadratic nonlinear terms with respect  to   the semi-groups associated to the linear part of the system    \eqref{SHGSA}.

First, recall the work of Bourgain in \cite{BOU98} where a refinement of Strichartz estimates in the Euclidean case was derived. Considering $ u_0 $, $ v_0 $  localized in frequency on the sets 
 $\{ \xi \in \R^{d} : |\xi| \sim N_1\}$ and 
 $\{ \xi \in \R^{d} : |\xi| \lesssim  N_2\}$ respectively, with
$N_2 \leq N_1$, the author in \cite{BOU98} proved the following estimate
\begin{equation}\label{BiliBou}
\| e^{i t \Delta} u_0  e^{i t \Delta} v_0\|_{L^{2} (\R \times \R^{d})} \lesssim \frac{N_2^{\frac{d-1}{2}}}{N_1^{\frac{1}{2}}}
\| u_0\|_{L^{2} ( \R^{d})} \| v_0\|_{L^{2} ( \R^{d})} \lesssim \min(N_1, N_2)^{\frac{d}{2} - 1} \| u_0\|_{L^{2} ( \R^{d})} \| v_0\|_{L^{2} ( \R^{d})}.
\end{equation}
Notice that if  $u_0 = v_0$, $N_1 = N_2$ and  $d = 2$, we can use the usual Littlewood-Paley decomposition 
on  $\R^{d}$ to obtain the already known Strichartz estimate
\begin{equation}
\left(  \int_{\mathbb{R}} \left( \int_{\mathbb{R}^{2}} |e^{i t \Delta } u_0(x)|^{4} dx \right) dt\right)^{1/4} \lesssim \|u_0\|_{L^{2}(\mathbb{R}^{2})}.
\end{equation}

As we are working on a compact Zoll manifolds considering two different groups, we need to derive an analogue of \eqref{BiliBou} that fits in our context.   Before entering to the details, we introduce the following  more general definition and list the known results in the literature.

\begin{defini}
 Let $P$, $Q$ be differential operators on   $M$  of order  $n$, $m$ respectively. 
We say that the associated semi-groups   $e^{i t P}, e^{i t Q}$ satisfy a
\textit{ bilinear estimate } of order  $s_0 := s_0(M) \geq 0$ if, for any $u_0$ localized in a frequency
 $N_1$ and  $v_0$ localized on a frequency  $N_2$ one has that
\begin{equation}\label{BiliMar}
\| e^{i t P} u_0  e^{i t Q} v_0\|_{L^{2} (I \times M)} \leq C(I) \min(N_1, N_2)^{s_0} 
\| u_0\|_{L^{2} (M)} \| v_0\|_{L^{2}(M)},
\end{equation}
where $I \subset \R$  is a finite interval. 
\end{defini}

In what follows we  list some works on the bilinear estimates \eqref{BiliMar} obtained in the context of compact Riemannian manifolds. 

\begin{center}
    \begin{tabular}{| l| l| l| l| l|l|l|l| l|p{3cm} |}
    \hline
   {\bf  Compact manifold $(M^{d},g)$} &{\bf Group 1}& {\bf Group 2} &  {\bf $s_0(M)$} & {\bf Authors}\\ \hline
    $ \mathbb{S}^{2}$ & $e^{i t \Delta}$& $e^{i t \Delta}$ & $\frac{1}{4}+$ & Burq, G\'erard, Tzvetkov \cite{BGT2}\\ \hline
     $ \mathbb{S}^{3}$  & $e^{i t \Delta}$ & $e^{i t \Delta}$ & $\frac{1}{2}+$ &Burq, G\'erard, Tzvetkov \cite{BGT1} \\ \hline
    $\mathbb{S}^{2} \times \mathbb{S}^{1} $ & $e^{i t (\Delta + \partial_{xx})}$ & $e^{i t (\Delta + \partial_{xx})}$ & $\frac{3}{4}+$ &Burq, G\'erard, Tzvetkov \cite{BGT1}  \\ \hline
     $\mathbb{T}^{2}, \mathbb{T}^{3} $  & $e^{i t \Delta}$ & $e^{i t \Delta}$ & $0+, \frac{1}{2}+$ &Bourgain \cite{BOU93A, BOU93B,BOU2007} \\ \hline
      $\tilde{\mathbb{T}}^{3} = \mathbb{R}^{3}/\Pi_{j = 1}^{3} (a_j \mathbb{Z})$ ($a_j \in \mathbb{R} \setminus \mathbb{Q}$)  & $e^{i t \Delta}$ & $e^{i t \Delta}$ & $\frac{2}{3}+$ & Bourgain \cite{BOU2007}\\ \hline
      $ \partial M \neq \emptyset$ ($d = 2$)
       & $e^{i t \Delta}$ & $e^{i t \Delta}$ & $\frac{2}{3}+$ &  Blair, Smith, Sogge \cite{BSS2012,JIA2011} \\ \hline
     Compact manifolds of $d \geq 2$   & $e^{i t \Delta}$& $e^{\pm i t |\nabla|} $ & $\frac{d - 1}{2}$ & Hani \cite{HAN2012} \\ \hline   
     $\tilde{\mathbb{T}}^{d} = \mathbb{R}^{d}/\Pi_{j = 1}^{d} (a_j \mathbb{Z})$ ($a_j \in \mathbb{R} \setminus \mathbb{Q}$) & $e^{i t \Delta}$& $e^{i t \Delta}$ & $\frac{d-2}{2}+$ & Fan, Staffilani, Wang, Wilson \cite{FSWW2018}\\ \hline
    \end{tabular}
\end{center}
\bigskip

With the above information in mind, as in the case of  the single NLS equation,  we plan to obtain  the bilinear estimates for the (now mixed) semi-groups
\begin{equation}\label{UVGroup}
(V(t), U_{\sigma}(t)) = \Big( e^{i t (\Delta - 1)}v_{0}, e^{i t (\frac{1}{\sigma}\Delta - \frac{\alpha}{\sigma})}u_{0}\Big),
\end{equation}
that describe the solution to the linear part
\begin{equation}\label{LSMod}
\begin{cases}
i \partial_{t} v + \Delta v - v = 0,\\
i \partial_{t} u + \frac{1}{\sigma}\Delta u - \frac{\alpha}{\sigma}u= 0, \\
(v(0), u(0)) =(v_0,u_0),
\end{cases}
\end{equation}
 associated to  the IVP \eqref{SHGSA}.
 
To better investigate the properties of the operators involved in the expressions of these semigroups, let us consider the operator
\[
- \Delta_{\delta} := -\delta \Delta + \gamma I = -\delta \Delta + \gamma .
\]
If  $(e_k, \mu_k)$ is an eigenpair of $-\Delta$, it is easy to check that $(e_k, \delta \mu_k + \gamma)$ is an eigenpair of  $-\Delta_{\delta}$. 
%In fact, 
%\[
%-\Delta_{\delta} e_k = \delta ( -\Delta e_k) + \gamma e_k = (\delta \mu_k + \gamma)e_k.
%\]
%On the other hand, if $\mu$ is an eigenvalue of $- \Delta_{\delta}$, there exists $u \neq 0$ such that
%\[
%- \Delta_{\delta} u = \mu u. 
%\]
%As $u \in L^{2}$, we can write $u = \sum_{k = 0}^{\infty} c_{k} e_{k}$ (expansion with respect to the eigenfunctions of $- \Delta$). Using $L^2-$ orthogonality,  and comparing the coefficients we can conclude that  $\mu = \delta \mu_{k}  + \gamma $ for some $k \in \mathbb{N}$. 
Considering  the case where $\delta = \frac{1}{\sigma}$ and $\gamma = -\frac{\alpha}{\sigma} $, we have 
\[
 e^{i t (\frac{1}{\sigma}\Delta - \frac{\alpha}{\sigma})}u_{0}   = e^{-i t\frac{\alpha}{\sigma}} e^{i t \frac{1}{\sigma}\Delta} u_{0}.
\]
Thus, considering  the product in $L^2$ with the semigroup generated by the choice of the parameters 
 $\beta = 1$ and $\gamma = -1$,  we have 
\begin{equation}\label{RenormUV}
\Big\| e^{i t (\frac{1}{\sigma}\Delta - \frac{\alpha}{\sigma})}u_{0} \mbox{ } e^{i t (\Delta - 1)}v_{0} \Big\|_{L^{2}(I \times M)} = \Big\| e^{i t \frac{1}{\sigma}\Delta}u_{0} \mbox{ } e^{i t \Delta}v_{0} \Big\|_{L^{2}(I \times M)}. 
\end{equation}
This shows that the product of the mixed groups is not affected by the  factors 
$e^{- it \frac{\alpha}{\sigma}}, \; e^{-it}$, where $I \subset \mathbb{R}$.

\subsection{Bilinear estimates on compact manifolds }

We start with the  following result that provides a bilinear estimate for localized functions on different regimes of frequency, in the case of compact manifolds.

\begin{prop}(\cite{BGT1} p. 261). \label{BiliSurf}
Let $M$ be a compact smooth manifold without boundary of dimension $d$. Let  $\chi \in C_0^{\infty}(\mathbb{R})$ and $s_0(d)$ defined in  \eqref{BiliSdb} . Given
$A > 0$, introduce the following approximated spectral projector, via functional calculus:
\[
\chi_A = \chi(\sqrt{- \Delta} - A). 
\]
Then, there is $C > 0$ such that, for all $\lambda, \mu \geq 1$, and any $f, g\in L^2(M)$, 
\begin{equation}\label{BiliSurfA}
\| \chi_\lambda f \chi_\mu g \|_{L^2(M)} \leq C \min(\lambda, \mu)^{s_0(d)} \|f\|_{L^2(M)} \|g\|_{L^2(M)}.
\end{equation}
\end{prop}

Next, we  use Proposition \ref{BiliSurf} to  characterize the localization operator $\chi_A$. 
\begin{enumerate}
\item[$\bullet$]  Choose  $\chi \in C_0^{\infty}(\mathbb{R})$ be such that $\supp \chi \subset(0, 2 \delta_0)$, 
$\chi \geq 0$ and  $\chi = 1$ in $[\delta_0/2, 3 \delta_0/2]$. Take $\sqrt{\mu_{k_0}} \in \Spec(\sqrt{- \Delta})$. As the spectrum of 
$\sqrt{- \Delta}$ is discrete, there exists $\delta_0 > 0$ such that
\[
\Spec(\sqrt{- \Delta}) \cap \left([\sqrt{\mu_{k_0}} - \delta_0, \sqrt{\mu_{k_0}} + \delta_0]=:I_0\right) \neq \emptyset.
\]
 Let  $A_0  = \sqrt{\mu_{k_0}} - \delta_0$, 
then $\chi( t - A_0) > 0$ if  $t \in I_0$. Moreover,  $\chi(t - A_0) = 1$ if
\[
 t \in [\sqrt{\mu_{k_0}} - \delta_0/2, \sqrt{\mu_{k_0}} + \delta_0/2] =: I'_0 \subset I_0. 
\]
Thus,
\[
\chi(\sqrt{- \Delta} - A_0) f  = \sum_{k \in \mathbb{N}} \chi(\sqrt{\mu_k} - (\sqrt{\mu_{k_0}} - \delta_0)) P_k f = P_{k_0} f. 
\]

Choose $\lambda = A_0$ and  $\mu = \tilde{A_0} = (\sqrt{\mu_{\tilde{k}_0}} - \tilde{\delta}_0) $,  with
$f = P_{k_0}$ and $g  = P_{\tilde{k}_0}$. Then from \eqref{BiliSurfA} we have
\begin{equation}\label{BiliSurfB}
\| P_{k_0} f P_{\tilde{k}_0} g \|_{L^2(M)} \leq  \min(\sqrt{\mu_{k_0}}, \sqrt{\mu_{ \tilde{k}_0}})^{\frac{1}{4}} 
\| P_{k_0} f \|_{L^2(M)} \| P_{\tilde{k}_0} g \|_{L^2(M)}.
\end{equation}

\item[$\bullet$] Now, let $\chi \in C_0^{\infty}(\mathbb{R})$ with $\supp \chi \subset(0, B)$ ($B > 0$), 
$\chi \geq 0$ and  $\chi = 1$ in $[\delta, B - \delta]$ for some $0 < \delta < B/2$. Then, for $A > 0$
\[
\chi(\sqrt{- \Delta} - A )f  = \sum_{k = 0 }^{+ \infty} \chi(\sqrt{\mu_k} - A ) P_k f , 
\]
where $\chi(\sqrt{\mu_k} - A ) > 0$ only for the values of   $k \in \mathbb{N}$ such that $ \sqrt{\mu_k} \in (A, A + B)$ and  $\chi(\sqrt{\mu_k} - A ) = 1$ 
if $\sqrt{\mu_k} \in [ A + \delta,  A + B - \delta]$. Therefore, we can write
\begin{equation}\label{SpecLoc}
\chi_{A}f  = \sum_{ k: \sqrt{\mu_k} \in [ A + \delta,  A + B - \delta] } P_k f  +  \sum_{k : \sqrt{\mu_k} \in (A, A + \delta) \cup (B + A - \delta, A+ B) } \chi(\sqrt{\mu_k} - A ) P_k f.
\end{equation}

Observe that $\Spec(\sqrt{- \Delta}) \cap (A, A + B) = \{ \sqrt{\mu_{k_{\ell}}} \}_{\ell = 1}^{n} $, where $ \sqrt{\mu_{k_{\ell}}} \leq  \sqrt{\mu_{k_{\ell + 1}}} $ for $1 \leq \ell \leq n - 1$. Therefore, if 
$ 0 < \delta < \min \{ | \sqrt{\mu_{k_1}} - A|,  | \sqrt{\mu_{k_n}} - (A + B)|\}$,  the second sum on the right hand side of \eqref{SpecLoc} vanishes. Consequently, the spectral projectors for which \eqref{BiliSurfA} is true are of the form 
\begin{equation}\label{BiliSurfC}
\chi_{A}f  = \sum_{ k: \sqrt{\mu_k} \in [ A + \delta,  A + B - \delta] } P_k f  = \sum_{\ell = 1}^{n} P_{k_{\ell}} f. 
\end{equation}
\end{enumerate}

\subsection{ Bilinear estimates in spheres}
We begin this section with a brief revision on the concept of Spherical Harmonics. For a complete introduction on the subject, see \cite{AH2012, DX2013}.  Let $P = P(x_{0}, x_{1}, ..., x_{d})$ be a homogeneous polynomial of degree  $k$ in $\mathbb{R}^{d + 1}$. If  $P$
is  harmonic, that is, $\sum_{i} \frac{\partial^{2} P}{\partial x_{i}^{2}} = 0$ then it can be shown that the restriction $ H := P \mid_{\mathbb{S}^{d}}$ 
satisfies $ \Delta_{\mathbb{S}^{d}} H  + \mu_{k} H = 0$, where
\begin{equation}\label{SpecSd}
 \mu_{k} := k (k + d - 1)  \qquad k \in \mathbb{N}, 
\end{equation} 
 denote the eigenvalues of $-\Delta_{\mathbb{S}^{d}}$. These functions are called spherical harmonics of degree $k$, which will be denoted by $H_{k}$. 

Next, we state a result that provides bilinear spectral estimates in  $d$-dimensional spheres,  $d \geq 2$.

\begin{prop}\label{BiliSd} (\cite{BGT1, BGT2}).
Let  $H_{k}, \widetilde{H}_{\ell}$ be spherical harmonics of degrees  $k, \ell \geq 1$ on $\mathbb{S}^{d}$ ($d \geq 2$). Then, 
\begin{equation}\label{BiliSda}
\| H_{k} \widetilde{H}_{\ell}\|_{L^{2}(\mathbb{S}^{d})} \leq C (\min (k, \ell))^{s_{0}(d)} \|H_{k}\|_{L^{2}(\mathbb{S}^{d})} \| \widetilde{H}_{\ell}\|_{L^{2}(\mathbb{S}^{d})},
\end{equation}
where $s_0(d)$ is defined in  \eqref{BiliSdb}. 
\end{prop}

\subsection{From the spectral estimates to the evolution estimates}

In this subsection, we will use the spectral bilinear estimates     \eqref{BiliSurfA} and \eqref{BiliSda} to infer
evolution bilinear estimates  involving
the interaction of the semigroups given in \eqref{RenormUV} by
\begin{equation}\label{SemiGPS}
S_{1/\sigma}(\pm t) := e^{ \pm i \frac{t}{\sigma} \Delta}  \qquad\mbox{ and }\qquad S(t) := e^{i t \Delta}.
\end{equation}

%%%%%%%%%%%%%%%%%%%%%%%%%%%%%%%%%%%%%%%%%%%%%%%%%%%%%%%%%%%%%%
We emphasize here that the choice of the  $\pm$ signs in \eqref{SemiGPS} is due to the presence of complex conjugate in the expressions we are going to deal with and has no relation to the choice of  $(+, +)$ signs that we made in the introduction. Notice  that, on the sphere $ (\mathbb{S}^d, can) $ we can take advantage of the precise knowledge
of the spectrum of the operator  $-\Delta_{\mathbb{S}^{d}}$ (see \eqref{SpecSd}). After obtaining the  bilinear estimates 
in the case of the spheres, we will extend those results to the case of Zoll manifolds.
Since the localization of the eigenvalues of the operator $-\Delta_{Zoll}$ is also well understood
(see Proposition \ref{ZollSpec} below), we will extend the estimate in this case by introducing a suitable abstract perturbation of the Laplacian to reduce to an analysis that is already carried out  in the case of the sphere $\mathbb{S}^{d}$.  %We will also derive 
%a generalization of the bilinear estimates for the case of the spheres $\mathbb{S}^d $,  $ d \geq 3 $,
%taking advantage  of the precise knowledge of the  spectrum.

The following lemma will be important in our argument which is proved using the estimate for the number of divisors of a natural number. 
%%%%%%%%%%%%%%%%%%%%%%%%%%%%%%%%%%%%%%%%%%%%%%%%%%%%%%%%%%%%%
\begin{lema}\label{NTLema} 
For every  $\varepsilon > 0$, there is  $C_{\varepsilon} > 0$
such that, given  $m \in \mathbb{Z} \setminus\{0\}$ and a positive integer $N$, 
\[
\sharp \{(x, y) \in \mathbb{N}^{2} \mid N \leq x \leq 2 N,\;\; x^{2} \pm y^{2} = m \} \leq C_{\varepsilon} N^{\varepsilon}. 
\]
\end{lema}

\begin{proof}
The proof for the case  $x^2 + y^2 = m$ is given in \cite{BGT2} p. 207. Thus, it remains to prove the hyperbolic case $x^2 - y^2 = m$. 

First, we consider the case when $m > 0$.  Since $ N \leq x \leq 2 N$, we have
$$
N^2 - m \leq x^2 - m \leq 4 N^2 - m,
$$
and consequently $y \in [0, 2N]$. Let $\tilde{x} = x - y$ and $\tilde{y} = x + y$. Note that $\tilde{x} \tilde{y} = m$ is possible  only if $\tilde{x} > 0$. Hence, it  suffices  to estimate
\begin{equation}\label{SetA}
 A_{m} := \sharp \{ (\tilde{x}, \tilde{y}) \mid  0 \leq \tilde{x} \leq 2N, \mbox{ }  N \leq \tilde{y} \leq 4 N \;\;  \mbox{and} \;\;  \tilde{x} \tilde{y} = m \}.    
\end{equation}
The number of divisors of a natural number $n$, denoted by $d(n)$, satisfies  $d(n) = O(n^{\epsilon})$. Thus,  one obtains  $A_{m} \lesssim d(8 N^2) = O(N^{\varepsilon})$. 

Next, we consider the case when $m < 0$. In this case, we write
\[
y^2 = x^2 - m = x^2 + m' , 
\]
where $m' := -m >0$. Then, $\sqrt{N^2 + m'} \leq y \leq \sqrt{4 N^2 + m'}$ and so,
$y\in I_{m'}: = \left[\sqrt{N^2 + m'} , \sqrt{4 N^2 + m'}\right]$. 

If $m' < 3 N^4$ then $y \in [N, 4N^2]$ and we may reduce to the same analysis of $A_m$. 

If $m' \geq 3 N^4$,  then we simply observe that $y$ 
takes values in  $I_{m'}$,  which  contains at most one integer,  since $|I_{m'}| < 1$. The proof of the lemma is completed
taking into account that if we fix $y$ then $x$ can not take more than one
value.
\end{proof}

Also, we will use the following result involving Fourier series.

\begin{lema}\label{BouLema}(\cite{BGT1} p. 289).
Let  $A \subset{\mathbb{R}}$ be a countable set. Then,  for every
$T > 0$ there is  $C_{T} > 0$ such that for every sequence  $(a_{\lambda})$ indexed by  $A$, one has 
\[
\left\|\sum_{\lambda \in A} a_{\lambda} e^{i \lambda t}  \right \|_{L^{2}(0, T)} \leq C_{T} \left(\sum_{m \in \mathbb{Z}} \left(\sum_{\lambda : |\lambda - m| \leq 1/2} |a_{\lambda}|\right)^{2} \right)^{1/2}. 
\]
\end{lema}

In what follows,  we  use the Lemmas \ref{NTLema} and \ref{BouLema} and the estimate \eqref{BiliSda} to derive  the bilinear Strichartz estimate on spheres.
%%%%%%%%%%%%%%%%%%%%%%%%%%%%%%%%%%%%%%%%%%%%%%%%%%%%%%%%%%%%%%%%%%
\begin{prop}\label{BiliEvo}(Bilinear Strichartz estimate on $\mathbb{S}^{d}$). Let $d \geq 2$ and
$s > s_0(d)$. 
 Consider the  semi-groups given by  \eqref{SemiGPS} with $\sigma = \frac{\beta}{\theta}$ where
$\beta, \theta \in \{n^{2}: n \in \mathbb{N} \}$. 
Then, there is   $C > 0$ such that for any  $u_{0}, v_{0} \in L^{2}(\mathbb{S}^{d})$ satisfying the spectral localization conditions 
 $u_{0}^{N} = u_{0}$, $v_{0}^{L} = v_{0}$,
one has 
 \begin{equation}\label{BiliEvoa}
    \|e^{ \pm i  \frac{t}{\sigma}\Delta}u_0 \mbox e^{i t \Delta} v_0\|_{L^{2}((0,1)_{t} \times \mathbb{S}^{d})} \leq C (\min(N, L))^{s} \|u_0\|_{L^{2}(\mathbb{S}^{d})}\|v_0\|_{L^{2}(\mathbb{S}^{d})}.
\end{equation}
\end{prop}
\noindent
\begin{prova} 
Initially consider $d = 2$. Let  us consider the  `` - "  sign,  the  proof  for the \;``+" sign  will follow in a similar way.

 Using the series expansion associated with the semi-groups 
 \eqref{SemiGPS} and the fact that  $u_{0}^{N} = u_{0}$, $v_{0}^{L} = v_{0}$, we can write
\begin{equation}\label{BiliProd}
e^{- i \frac{t}{\sigma }  \Delta }u_{0} e^{i t \Delta}v_{0}  = \sum_{k:N \leq \langle \mu_k\rangle^{1/2} < 2 N } \sum_{ \ell :L \leq \langle \mu_\ell\rangle^{1/2} < 2 L }  e^{i t  (\frac{\mu_k}{\sigma} - \mu_\ell)} P_{k}u_{0} P_{\ell}v_{0} =: S^{NL}_{\sigma}(t).
\end{equation}

Applying the  Fubini's theorem, we get  
\begin{equation}
\|S^{NL}_{\sigma}\|^{2}_{L^{2}((0, 1) \times \mathbb{S}^{2})}  = \int_{0}^{1} \|S^{NL}_{\sigma}(t)\|^{2}_{ L^{2}(\mathbb{S}^{2})}  dt 
 = \int_{\mathbb{S}^{2}}  \| S^{NL}_{\sigma}\|_{L^2(0,1)}^{2}  dS.  
\end{equation}
 Now, from  Lemma \ref{BouLema}, we obtain
\[
\| S^{NL}_{\sigma}\|_{L^2(0,1)}^{2} \leq C \sum_{m \in \mathbb{Z}} \left( \sum_{ (k, \ell): |m - (\frac{\mu_k}{\sigma} - \mu_\ell)| \leq \frac{1}{2} \atop N \leq \langle \mu_k \rangle^{1/2}\leq 2N, L \leq \langle \mu_\ell \rangle^{1/2}\leq 2L} |P_ku_0 P_{\ell} v_0|\right)^{2}.
\] 
 Let
\begin{equation}\label{XX0}
\Lambda^{NL}(m) : = \left\{ (k, \ell) \in \mathbb{N}^{2}: \left|m - (\frac{\mu_k}{\sigma} - \mu_\ell)\right| \leq \frac{1}{2}; N \leq \langle \mu_k \rangle^{1/2} < 2N; L \leq \langle \mu_\ell \rangle^{1/2} < 2L \right\}.
\end{equation}
Applying the Cauchy-Schwarz inequality in the sum of $(k, \ell)$ and the 
triangle inequality for the $L^2(\mathbb{S}^2)$-norm, we obtain
\[
\|S^{NL}_{\sigma}\|^{2}_{L^{2}((0, 1) \times \mathbb{S}^{2})} \leq  C \sup_{m \in \mathbb{Z}} \sharp \Lambda^{NL}(m) 
\sum_{m \in \mathbb{Z}} \left( \sum_{(k, \ell) \in  \Lambda^{NL}(m) } \| P_{k}u_0 P_{\ell}v_0 \|_{L^{2}(\mathbb{S}^2)}^{2} \right).
\]
By Proposition \ref{BiliSd}, it follows that
\begin{equation}\label{XX1}
\begin{split}
\|S^{NL}_{\sigma}\|^{2}_{L^{2}((0, 1) \times \mathbb{S}^{2})} &\leq  C \min(N, L)^{1/2} \sup_{m \in \mathbb{Z}} \sharp \Lambda^{NL}(m) \\
& \times \sum_{m \in \mathbb{Z}} \left( \sum_{(k, \ell) \in  \Lambda^{NL}(m) } \| P_{k}u_0\|_{L^{2}(\mathbb{S}^2)}^{2} \| P_{\ell}v_0\|_{L^{2}(\mathbb{S}^2)}^{2} \right).
\end{split}
\end{equation}

On the one hand, considering  $\sigma := \frac{\beta}{\theta}$ with $\beta, \theta \in \{ n^2: n \in \mathbb{N} \}$ and   $\mu_j = j^2 + j$ ($j \in \mathbb{N}$), we can write 
\begin{equation}\label{NewSet}
\Lambda^{NL}(m) \subset \left\{(k, \ell) \in \mathbb{N}^{2}: \left| k_{\theta}^{2} - \ell_{\beta}^{2} - (4 \beta m + \theta - \beta) \right| \leq 2 \beta;  \frac{N}{2} \leq k \leq 2 N; \frac{L}{2} \leq \ell \leq 2 L \right \},
\end{equation}
 where $k_{\theta} := \sqrt{\theta} (2 k + 1)$ and  $\ell_{\beta} := \sqrt{\beta} (2 \ell + 1)$. Notice that
 the condition on $ \sigma $ is of technical character, since we must complete the square to obtain the expression of the set
 given in \eqref{NewSet} and also because we need to obtain an upper bound for its cardinality using Lemma \ref{NTLema}. Furthermore, we have
\begin{equation}\label{XX2}
\sharp \Lambda^{NL}(m) \leq \sharp \tilde{\Lambda}^{NL}(m)
\end{equation}
where
\[
\tilde{\Lambda}^{NL}(m) := \left\{(\tilde{k}, \tilde{\ell}) \in \mathbb{N}^{2}: \left| \tilde{k}^{2} - \tilde{\ell}^{2} - (4 \beta m + \theta - \beta) \right| \leq 2 \beta;  N_{\theta} \leq \tilde{k} \leq 4 N_{\theta}; L_{\beta} \leq \tilde{\ell} \leq 4 L_{\beta} \right \},
\] 
 $N_{\theta }:= \sqrt{\theta} (N + 1)$ and  $L_{\beta }:= \sqrt{\beta} (L + 1)$. 
 
  Let $\beta' \in ([-2 \beta, 2 \beta] \cap \mathbb{Z})$ then, using Lemma \ref{NTLema} uniformly with respect to    $m$, we obtain that the number of the pairs  $(\tilde{k}, \tilde{\ell})$, satisfying
\begin{equation}\label{N-X1}
 \tilde{k}^{2} - \tilde{\ell}^{2} - (4 \beta m + \theta - \beta) = \beta',
\end{equation}
is bounded by $O_{\epsilon}(\min(N, L)^{\epsilon})$. Since there are
 $4 \beta + 1$ possibilities for  $\beta'$, we have
\begin{equation}\label{XX3}
\sharp \tilde{\Lambda}^{NL}(m) \leq (4 \beta + 1) O_{\epsilon}(\min(N, L)^{\epsilon}).
\end{equation}
From  \eqref{XX1}, \eqref{XX2} and  \eqref{XX3} it follows that
\begin{equation}\label{XX4}
\|S^{NL}_{\sigma}\|^{2}_{L^{2}((0, 1) \times \mathbb{S}^{2})} \lesssim_{\epsilon, \beta}   \min(N, L)^{1/2 + \epsilon} 
\sum_{m \in \mathbb{Z}} \left( \sum_{(k, \ell) \in  \Lambda^{NL}(m) } \| P_{k}u_0\|_{L^{2}(\mathbb{S}^2)}^{2} \| P_{\ell}v_0 \|_{L^{2}(\mathbb{S}^2)}^{2} \right).
\end{equation}

On the other hand, according to \eqref{XX0}, we have 
 $\Lambda^{NL}(m) = \emptyset$ if 
\[
 m \notin \Big[- 4 L^2 - 1, \frac{4N^2}{\sigma} + 1\Big] \cap \mathbb{Z}. 
\]
Thus, to bound the summations in the RHS of \eqref{XX4}
we need to consider 
\[
\sum_{(k, \ell) \in \Big( \bigcup_{ m \in \big[- 4 L^2 - 1,  \frac{4N^2}{\sigma} + 1\big] \cap \mathbb{Z}} \Lambda^{NL}(m)\Big) } \| P_{k}u_0\|_{L^{2}(\mathbb{S}^2)}^{2} \| P_{\ell}v_0 \|_{L^{2}(\mathbb{S}^2)}^{2}.
\]
But, 
\[
 \bigcup_{ m \in \big[- 4 L^2 - 1, \frac{4N^2}{\sigma}  + 1\big] \cap \mathbb{Z}} \Lambda^{NL}(m) \subset 
 \left\{ (k, \ell) \in \mathbb{N}^{2}: N \leq \langle \mu_k \rangle^{1/2} < 2N; L \leq \langle \mu_\ell \rangle^{1/2} < 2L \right\}.
\]
Therefore, the summations on the right side of \eqref {XX4} can be bounded by
\[
\left(\sum_{k : N \leq \langle \mu_k \rangle^{1/2}\leq 2N}\| P_{k}u_0\|_{L^{2}(\mathbb{S}^2)}^{2} \right) 
\left(\sum_{\ell : L \leq \langle \mu_\ell \rangle^{1/2}\leq 2L}\| P_{\ell}v_0\|_{L^{2}(\mathbb{S}^2)}^{2} \right) =
 \| u_0^{N}\|_{L^{2}(\mathbb{S}^2)}^{2} \| v_0^{L}\|_{L^{2}(\mathbb{S}^2)}^{2}. 
\]
That is,
\begin{equation*}
\|S^{NL}_{\sigma}\|^{2}_{L^{2}((0, 1) \times \mathbb{S}^{2})} \lesssim_{\epsilon, \beta}   \min(N, L)^{1/2 + \epsilon} 
\| u_0^{N}\|_{L^{2}(\mathbb{S}^2)}^{2} \| v_0^{L}\|_{L^{2}(\mathbb{S}^2)}^{2}.
\end{equation*}
By taking $\epsilon: = 2s - 1/2$ we finish the proof for $d=2$. 

For $d \geq 3$,  we follow the same lines of the proof given for  $d= 2$.  First, exploiting the knowledge of the spectrum of  $- \Delta_{\mathbb{S}^d}$ (see \eqref{SpecSd}),  we need to analyse the cardinality of the set 
\begin{equation}\label{XX0-d}
\Lambda^{NL}_d(m) : = \Big\{ (k, \ell) \in \mathbb{N}^{2}: \Big|m - (\frac{\mu_k}{\sigma} - \mu_\ell)\Big| \leq \frac{1}{2}; N \leq \langle \mu_k \rangle^{1/2}\leq 2N; L \leq \langle \mu_\ell \rangle^{1/2}\leq 2L \Big\},
\end{equation}
with $\mu_k=k(k+d-1)$ and $\mu_{\ell}=\ell(\ell+d-1)$.

Analogously to what was done previously, we must find a bound for the number of solutions of the inequality
\begin{equation}\label{DiophEqbb}
\Big|(\sqrt{\theta} (2k + d - 1))^{2} - (\sqrt{\beta} (2\ell + d - 1))^{2} - 4 \beta m + (\theta - \beta) (d -1)^2\Big| \leq 2 \beta,  
\end{equation}
where $k \in [N/2, N]$ and $\ell \in [L/2, L]$. 
Using Lemma \ref{NTLema} and the spectral bilinear estimate given by Lemma \ref{BiliSd} for  $d \geq 3$, we obtain
 \begin{equation*}\label{BiliEvob}
    \|e^{ \pm i  \frac{t}{\sigma}\Delta}u_{0} \mbox e^{i t \Delta} v_{0}\|_{L^{2}((0,1)_{t} \times \mathbb{S}^{d})} \leq C (\min(N, L))^{s} \|u_{0}\|_{L^{2}(\mathbb{S}^{d})}\|v_{0}\|_{L^{2}(\mathbb{S}^{d})}.
\end{equation*}
for all $s > s_0(d)$ where $s_0(d)$ is defined in  \eqref{BiliSdb}.
\end{prova}

%%%%%%%%%%%%%%%%%%%%%%%%%%%%%%%%%%%%%%%%%%%%%%%%%%%%%%%%%%%%%%%%%%%%%%%%

The following lemma provides a reformulation of the bilinear 
Strichartz estimates in terms of bilinear estimates involving 
the mixed spaces $X^{s,b}_{\delta}$ and $X^{s, b}$.

\begin{lema}\label{Eq00}
Let  $s > 0$ and  $\sigma > 0$ be as in the statement of the Proposition  \ref{BiliEvo}. The following statements are equivalent
\begin{enumerate}
\item[$(i)$] For any $u_0, v_0 \in L^{2}(M)$ satisfying  
$u_{0}^{N} = u_0$ and  $v_{0}^{L} = v_0$, one has  
\begin{equation}\label{Eq00a}
    \|S_{1/\sigma}(\pm t)u_{0} \mbox{ } S(t) v_{0}\|_{L^{2}((0,1)_{t} \times M)} \leq C (\min(N, L))^{s} \|u_{0}\|_{L^{2}(M)}\|v_{0}\|_{L^{2}(M)}.
\end{equation}
\item[$(ii)$]
 For all $b > 1/2$ and any functions $u_{1} \in X_{1/\sigma}^{0, b}(\mathbb{R} \times M)$ and  $u_{2} \in X^{0, b}(\mathbb{R} \times M)$ satisfying $u_{1}^{N} =  u_{1}$ and  $u_{2}^{L} = u_{2}$, one has  
\begin{equation}\label{Eq00b}
\|\tilde{u_{1}} u_{2}\|_{L^{2}(\mathbb{R} \times M)} \leq C \min(N, L)^{s} \|u_{1}\|_{X^{0, b}_{1/\sigma}(\mathbb{R} \times M)} 
 \|u_{2}\|_{X^{0, b}(\mathbb{R} \times M)}, 
\end{equation}
where  $\tilde{u_{1}} = \overline{u_{1}}$ or  $ \tilde{u_{1}} = u_{1}$.
\end{enumerate} 
\end{lema}
\noindent
\begin{prova}
The proof follows with simple modification of Lemma 2.3 in \cite{BGT2} using the semigroups $S_{1/\sigma}(\pm t)$, $S(t)$ and  properties of the spaces $X^{s,b}_{1/\sigma}$, see also Proposition 4.3 in \cite{BGT1}. 
\end{prova}

Now, we prove that the bilinear Strichartz estimate stated in Proposition \ref{BiliEvo} can  be extended to the class of the $d$-dimensional  Zoll  manifolds. For this purpose, as we have already said, the following
result due to Colin de Verdi\`ere  \cite{CV1979} and Guillemin \cite{GUI1977} on the localization of the eigenvalues of the Laplacian on these manifolds (which we denote by $- \Delta_{Zoll}$ just to make it more explicit) plays  crucial role. 
It is worth mentioning that this localization property holds in higher dimensions and not only in the case of surfaces $(d = 2)$ (see for example \cite{WEI1997,ZEL1996,ZEL1997}). 

\begin{prop}\label{ZollSpec}
If the geodesics of  $M$ are $2 \pi$-periodic\footnote{This renormalization is done by applying a dilation in the Riemannian metric $ g $, that is, given an appropriate scalar $ \lambda \in (0, + \infty) $ we consider the metric $\lambda g$ instead of $g$. }, there are $Z_0 \in \mathbb{N}$ and  $ E > 0$
 such that the spectrum of $- \Delta_{Zoll} $ is contained in  $ \bigcup_{ k = 1}^{+ \infty} I_k$, where
\begin{equation*}\label{Ik}
I_k : = \Big[\Big(k + \frac{Z_0}{4}\Big)^2 - E, \left(k + \frac{Z_0}{4}\right)^2 + E \Big].
\end{equation*}
\end{prop}

\begin{prop}\label{BiliEvo2}(Bilinear evolution estimates on Zoll manifolds). Let $M$ be a $d$-dimensional Zoll manifold  and
$s > s_0(d)$  be given. 
 Consider the  semi-groups given by \eqref{SemiGPS} with $\sigma = \frac{\beta}{\theta}$ where
$ \beta, \theta \in \{n^{2}: n \in \mathbb{N} \}$. 
Then, there is  $C > 0$ such that for any  $u_{0}, v_{0} \in L^{2}(M)$ satisfying
the  spectral localization conditions $u_{0}^{N} = u_{0}$, $v_{0}^{L} = v_{0}$,
one has 
\begin{equation}\label{BiliZoll}
    \|e^{ \pm i  \frac{t}{\sigma}\Delta}u_{0} \mbox e^{i t \Delta} v_{0}\|_{L^{2}((0,1)_{t} \times M)} \leq C (\min(N, L))^{s} \|u_{0}\|_{L^{2}(M)}\|v_{0}\|_{L^{2}(M)}.
\end{equation}
\end{prop}
\noindent
\begin{dem}
Denote by $(e_n)_n$ the sequence of eigenfunctions of $- \frac{1}{\sigma} \Delta_{Zoll}$ associated 
to the eigenvalues $\lambda_n^{\sigma}$ (counting its multiplicity). Using the Proposition \ref {ZollSpec}, it follows that
$\Spec(- \frac{1}{\sigma} \Delta_{Zoll}) \subset \bigcup_{ k = 1}^{+ \infty} I^{\sigma}_k $, where $I^{\sigma}_k : = \frac{1}{\sigma} I_k$. Observing the expressions of the $I_k$ intervals, it can be inferred  that they are not necessarily disjoint. In this way, it is easy to see that, if we take
\[
K_0 : = \min \Big\{  k \in \N \mid k > E - \frac{Z_0}{4} - \frac{1}{2} \Big\},
\]
we have  $I^{\sigma}_m \cap I^{\sigma}_\ell = \emptyset$, whenever $m \neq l$ and   $m, \ell \geq K_0$.

Now, set $N_0$ as being the smallest natural number such that
\[
\lambda_{N_0}^{\sigma} \geq \left(K_0+ \frac{Z_0}{4}\right)^2- E .
\]
In this case, since the sequence of eigenvalues is non-decreasing, we have, for all $ n \geq N_0 $, that $ \lambda_ {\sigma} $
belongs to only one interval of the form $ I_k^{\sigma}$ (with 
$k \geq K_0$). Take 
$$L_{0} = \max \{ N_0, K_0\}.$$

Define  an abstract perturbation of the Laplacian on $L^{2}(M)$ by 
\begin{equation*}\label{AbstractL}
 - \frac{1}{\sigma} \Delta_{Zoll}^{\bigstar} e_n =\begin{cases}
\lambda_n^{\sigma} e_n  &\mbox{ if } n \leq L_0 \\
\frac{1}{\sigma}\left(k + \frac{Z_0}{4}\right)^2 e_n  &\mbox{ if } n > L_0 \mbox{ and  } \lambda_n^{\sigma} \in I_k^{\sigma}. \\
\end{cases}
\end{equation*}

According to Lemma \ref{Eq00}, to prove \eqref{BiliZoll}, it  suffices  to prove that  \eqref{Eq00b} is valid for  $s > s_0(M) = 1/4$ and   $X^{0,b}_{- \frac{1}{\sigma}\Delta_{Zoll}}$, $X^{0,b}_{ - \Delta_{Zoll}}$. But, according to Remark \ref{XsbObs},
\begin{equation*}\label{XsbEqua}
\begin{cases}
\| \cdot \|_{X^{0,b}_{ - \frac{1}{\sigma}\Delta_{Zoll}}}  \simeq \| \cdot \|_{X^{0,b}_{ - \frac{1}{\sigma}\Delta_{Zoll}^{\bigstar}}} \\
\| \cdot \|_{X^{0,b}_{ - \Delta_{Zoll}}} \simeq \| \cdot \|_{ X^{0,b}_{ - \Delta_{Zoll}^{\bigstar}}}.
\end{cases}
\end{equation*}
So, using the fact that $ (ii) $ implies $ (i) $ in the Lemma 
\ref{Eq00}, it is sufficient to prove that
\eqref{Eq00a} is valid with $ s> 1/4 $ and 
\begin{equation*}
\begin{cases}
S_{1/\sigma}(\pm t) = e^{ \pm i t \frac{1}{\sigma}\Delta_{Zoll}^{\bigstar}} \\
S(t) = e^{  i t \Delta_{Zoll}^{\bigstar}}.
\end{cases}
\end{equation*}

Replacing $ P_ {k} $ in \eqref{BiliProd} by\footnote {Note that this is a type of spectral projector for which the estimate \eqref {BiliSurf} is valid, and the same is true if we consider $\Pi_{I_{k}^{\sigma}}$.} $\Pi_{I_{k}} = \sum_{ n \in \N :\lambda_n \in \Spec(- \Delta_{Zoll}^{\bigstar}) \cap I_{k}} P_n$  and the set in \eqref{XX0-d} by
\begin{equation*}\label{XX0-Zoll}
\Lambda^{NL}_{Zoll}(m) : = \left\{ (k, \ell) \in \mathbb{N}^{2}: \left|m - (\frac{\mu_k}{\sigma} - \mu_\ell)\right| \leq \frac{1}{2}; N \leq \langle \mu_k \rangle^{1/2} < 2N; L \leq \langle \mu_\ell \rangle^{1/2} < 2L \right\},
\end{equation*}
with $\mu_k =(k + Z_0/4)^2$ and $\mu_{\ell}= (\ell + Z_0/4)^2$ and 
proceeding analogously as in the case of $\mathbb{S}^d$, we need to find a bound for the number of solutions of the inequality
\begin{equation*}\label{DiophZoll}
\left|(\sqrt{\theta} (4 k + Z_0))^{2} -   (\sqrt{\beta} (4 \ell + Z_0))^{2}   - 16 \beta m\right| \leq 8 \beta.
\end{equation*}
Thus, the proof given in Proposition \ref{BiliEvo} can be applied  to the   pair $(\pm \frac{1}{\sigma} \Delta_{Zoll}^{\bigstar}, \Delta_{Zoll}^{\bigstar})$ as well.
\end{dem}

\begin{obs}
In view of the result obtained in the Proposition \ref{BiliEvo2},  the equivalence of the estimates given in the Lemma \ref{Eq00} can be used for the case of  the Zoll manifolds whenever $s>s_0(d)$.
\end{obs}

 %%%%%%%%%%%%%%%%%%%%%%%%%%%%%%%%%%%%%%%%%%%%%%%%%%%%%%%%%%%%%

\section{ Bilinear estimates for quadratic interactions }\label{SectionDC4}

In this section, we will use the estimates obtained in the Propositions
\ref{BiliEvo} and \ref{BiliEvo2} and the equivalence given by Lemma 
\ref{Eq00} to obtain the bilinear estimates involving the quadratic interactions of the system \eqref{SHGSA}. More precisely, we prove the following result.

\begin{prop}\label{Bili1}(Bilinear estimates). Let $M$ be a $d$-dimensional Zoll manifold, $ s > s_0(d)$
and  $\sigma = \frac{\beta}{\theta}$ with $\beta, \theta \in \{ n^{2} : n \in \N \}$.
Then there exist  $(b,b') \in \mathbb{R}^{2}$ satisfying
$0 < b' < 1/2 < b$, $b + b' < 1$  and   $C > 0$ such that
\begin{equation}\label{Bili1a}
\|u_{1} \overline{u_{2}}\|_{X^{s, -b'}} \leq C \|u_{1}\|_{X^{s, b}_{1/\sigma}} \|u_{2}\|_{X^{s, b}} ,\quad \mbox{ if }\;  s_0<s<1,
\end{equation} 
 and 
\begin{equation}\label{Bili1b}
\|v_{1} \overline{v_{2}}\|_{X^{s, -b'}} \leq C \|v_1\|_{X^{s, b}_{1/\sigma}} \|v_2\|_{X^{1, b}}, \quad \mbox{ if } \; s \geq 1. 
\end{equation}
Also, there exist   $ (b_1,b'_1) \in \mathbb{R}^{2} $ satisfying  
$0 < b'_1 < 1/2 < b_1$, $b_1 + b_1' < 1$ and 
  $C > 0$ such that 
\begin{equation}\label{Bili2a}
\|u_{1} u_{2}\|_{X_{1/\sigma}^{s, -b_1'}} \leq C \|u_{1}\|_{X^{s, b_1}} \|u_{2}\|_{X^{s, b_1}},\quad \mbox{ if  } \;s_0<s<1
\end{equation}
and 
\begin{equation}\label{Bili2b}
\|v_{1} v_{2}\|_{X_{1/\sigma}^{s, -b_1'}} \leq C \|v_{1}\|_{X^{s, b_1}} \|v_{2}\|_{X^{1, b_1}}, \quad\mbox{ if  }\;  s \geq 1. 
\end{equation}

\end{prop}

To prove  Proposition \ref{Bili1}, we need a series of basic results
that we derive  in the following subsection.   

\subsection{Auxiliary results}

We  begin by proving a result which provides an estimate for the product of three functions in $M$ localized on three frequency intervals $I_{k}$, such that one of the intervals is very dislocated in relation to the others.

As shown in \cite{HAN2012B} page 1203, in the case of
 $M = \mathbb{T}^d$ or $\mathbb{S}^d$, we have
\[
\int_{M} e_{k_0}  \tilde{e_{k_1}}  \tilde{e_{k_2}} \tilde{e_{k_3}} dg = 0,
\]
if $\mu_{k_0} > \mu_{k_1} + \mu_{k_2} + \mu_{k_3}$, where $ \tilde{e_{k_j}} = e_{k_j}$ or  $\overline{e_{k_j}}$. Therefore, in these cases the lemma that we will announce  below is not necessary. This kind of property is often used when decompositions of function $u$ over $\mathbb{S}^d$ with respect to the spherical harmonics
of degree $ k $ are considered,
that is, decompositions of the type
\[
u(x) = \sum_{k} H_{k}(x) \in L^{2}(\mathbb{S}^{d}),
\]
for details see \cite {PG10} page 145 in the case of $ \mathbb{S}^4 $ and \cite{LAU2010} page 819 in the case of
$\mathbb{S}^3$. To deal with the case of  Zoll manifolds (where the cancellation property may fail), it is necessary to use the following result.

\begin{lema}\label{EqStar}
Let $M$ be a compact $d$-dimensional  Riemannian manifold. If there exist $C > 1$ such that
$C (\mu_{k_{1}} + \mu_{k_{2}})  \leq \mu_{k_0},$
then for all $p > 0$ there exists $C_{p} > 0$ such that, for all $w_{j} \in L^{2}(M)$, $j = 0, 1, 2$,
\[
\left|\int_{M} P_{k_{0}} w_{0} P_{k_{1}} w_{1} P_{k_{2}} w_{2} dg\right| \leq C_{p} \mu_{k_{0}}^{-p} \prod_{j = 0}^{2} \|w_{j}\|_{L^{2}(M)}.
\]
\end{lema}
\noindent
\begin{dem} The proof consists in reducing the estimate of  Theorem 4.2  of \cite{HAN2012B} to the case of the product of three projections.  More precisely in  \cite{HAN2012B}, (see \cite{BGT2} p. 198 for a different proof in the  case $d = 2$), it has been proved that if $C (\mu_{k_1} + \mu_{k_2} +\mu_{k_3}) \leq \mu_{k_0}$,
then for all $ \mathbb{N} \ni p > 0$ there exists $C_{p} > 0$ such that, for any $w_{j} \in L^{2}(M)$, $j = 0, 1, 2,3$,
\begin{equation}\label{XP3}
\left|\int_{M} P_{k_{0}} w_{0} P_{k_{1}} w_{1} P_{k_{2}} w_{2}  P_{k_{3}} w_{3} dg\right| \leq C_{p} \mu_{k_{0}}^{-p} \prod_{j = 0}^{3} \|w_{j}\|_{L^{2}(M)}.
\end{equation}
Thus, if we take $k_{3} = 0$ in \eqref{XP3}, we have  $\mu_{k_{3}} = \mu_{0} =0$ (first eigenvalue of  $- \Delta$ with eigenfunction
$e_{0} = 1$). Now, if $w_{3} = e_{0} = 1 \in L^{2}(M)$, we have
\[
P_{k_{3}} w_{3} = e_{k_{3}} \int_{M} w_{3} e_{k_{3}} dg  = e_{0} \int_{M} 1 \cdot e_{0} dg = \vol_{g}(M).
\]
Thus,  if $C (\mu_{k_1} + \mu_{k_2})  \leq \mu_{k_{0}}$,  one obtains
\[
\left|\int_{M} P_{k_0} w_{0} P_{k_1} w_1 P_{k_2} w_2 dg \right| \leq \widetilde{C}_{p} \mu_{k_0}^{-p} \prod_{j = 0}^{2} \|w_{j}\|_{L^{2}(M)},
\]
as required.
\end{dem}

The estimate obtained in the Lemma \ref{EqStar} will be crucial to establish estimates for an integral involving the product
of three functions localized  in different   frequency regimes under 
the condition $ N_0 \geq C (N_1 + N_2)$.

\begin{obs}\label{XP6}
From now on, to simplify the notation, we will use 
 $\langle \mu_{k}\rangle \sim N_{j}^{2}$ for $j = 0, 1, 2$ to denote  that $N_{j} \leq \langle \mu_{k} \rangle^{\frac{1}{2}} < 2 N_{j}$ and  $\langle \tau + \widetilde{\mu_{k}}\rangle \sim L_{j}$ to  indicate that $L_{j} \leq  \langle \tau + \widetilde{\mu_{k}}\rangle < 2 L_{j}$, where  $\widetilde{\mu_k} := \frac{\mu_k}{\sigma}$.
\end{obs}

\begin{lema}\label{LemaFund} Let $M$ be a compact $d$-dimensional  Riemannian  manifold.
Consider the expression 
\begin{equation}\label{XP5}
I^{NL} := \left|\int_{\mathbb{R} \times M} u_{0}^{N_{0}L_{0}} \overline{u_{1}^{N_{1}L_{1}}} u_{2}^{N_{2}L_{2}}  \right|,
\end{equation}
where 
\begin{equation}\label{Eq12X}
u_{j}^{N_{j}L_{j}}(t) = \frac{1}{2\pi} \sum_{ k: N_{j} \leq \langle \mu_{k}\rangle^{1/2} < 2 N_{j}} \langle \mu_{k}\rangle^{ \tilde{s}/2}  \int_{L_{j}\leq \langle \tau + \tilde{\mu_{k}}\rangle < 2 L_{j}} \langle \tau + \tilde{\mu_{k}}\rangle^{- \tilde{b}} \widehat{P_{k} w_{j}}(\tau) e^{i t \tau} d\tau,
\end{equation}
with 
\begin{equation}\label{Eq12Xa}
\begin{cases}
  \tilde{s} =  s \; \quad \mbox{ and } \;  \tilde{b} = b'  &\mbox{ if } j = 0;  \\
  \tilde{s} = -s  \; \mbox{ and } \; \tilde{b} = b   &\mbox{ if } j = 1, 2. 
  \end{cases}
\end{equation}
% $\tilde{s} =  s$ if  $j = 0$,  $\tilde{s} = -s$ if $j = 1, 2$, $\tilde{b} = b$ if  $j = 1,2 $  and   $\tilde{b} = b'$ if  $j = 0$. 
Let  $s > 0$. Then, there exist $C_{1} > 0$ and  $p > 0$ such that if   $N_{0} \geq C (N_{1}+ N_{2})$, we have
\begin{equation}\label{LemaFunda}
 I^{NL} \leq C_{1}  \frac{N_{0}^{s + 6 - p}}{ N_{1}^{s}  N_{2}^{s}}
 \frac{L_{2}^{1/2}}{L_{0}^{b'}L_{1}^{b} L_{2}^{b}}    \prod_{j = 0}^{2} \|w_{j}\|_{L^{2}(\mathbb{R} \times M)},
\end{equation}
where $ s + 6 - p < 0 $. 
\end{lema}
\noindent
\begin{prova}  Inserting the terms given in \eqref{Eq12X} for $j = 0, 1, 2$ in \eqref{XP5},  we obtain
\[
I^{NL} \leq \frac{1}{(2 \pi)^{3}} \sum_{k_{0} : N_{0}^{2} \sim \langle\mu_{k_{0}}\rangle} \sum_{k_{1} : N_{1}^{2} \sim \langle\mu_{k_{1}}\rangle} \sum_{k_{2} : N_{2}^{2} \sim \langle\mu_{k_{2}}\rangle} \frac{\langle \mu_{k_{0}}\rangle^{s}}{\langle \mu_{k_{1}}\rangle^{s}  \langle \mu_{k_{2}} \rangle^{s}} \times
\]
\begin{equation}\label{Eq16A}
\times \frac{1}{L_{0}^{b'}L_{1}^{b} L_{2}^{b}}\left|\int_{\mathbb{R} \times M} \left( \int_{D} \tilde{P}_{k_{0}} \overline{\widehat{w_{0}}(\tau_{0})} P_{k_{1}}\widehat{w_{1}}(\tau_{1}) \tilde{P}_{k_{2}} \overline{\widehat{w_{2}}(\tau_{2})} e^{i t (  \tau_{1} - \tau_{2} - \tau_{0})} d \tau_{0} d \tau_{1} d \tau_{2}\right) dg dt \right|,
\end{equation}
where
\[
D := \{ (\tau_{0}, \tau_{1}, \tau_{2}) \in \mathbb{R}^{3} : \langle \tau_{0} + \mu_{k_{0}}\rangle \sim L_{0}, \langle \tau_{1} + \frac{\mu_{k_{1}}}{\sigma}\rangle \sim L_{1} , \langle \tau_{2} + \mu_{k_{2}}\rangle \sim L_{2} \}.
\]
%%%%%%%%%%%%%%%%%%%%%%%%%%%%%%%%%%%%%%%%%%%%%%%%%%%%%%%%%%%

Denote 
\[
I(\tau_0, \tau_1, \tau_2) := \int_{\mathbb{R} \times M} \left( \int_{D} \tilde{P}_{k_{0}} \overline{\widehat{w_{0}}(\tau_{0})} P_{k_{1}}\widehat{w_{1}}(\tau_{1}) \tilde{P}_{k_{2}} \overline{\widehat{w_{2}}(\tau_{2})} e^{i t (  \tau_{1} - \tau_{2} - \tau_{0})} d \tau_{0} d \tau_{1} d \tau_{2}\right) dg dt.
\]
 Let $\delta$ be the Dirac measure. As in \cite{BGT2} p. 200, we  can rewrite 
\[
I(\tau_0, \tau_1, \tau_2) =  \int_{D \times M} \tilde{P}_{k_{0}} \overline{\widehat{w_{0}}(\tau_{0})} P_{k_{1}}\widehat{w_{1}}(\tau_{1}) \tilde{P}_{k_{2}} \overline{\widehat{w_{2}}(\tau_{2})}  dg d \tau_{0} d \tau_{1} d \tau_{2}  \delta(- \tau_0 + \tau_1 - \tau_2) 
\]

%%%%%%%%%%%%%%%%%%%%%%%%%%%%%%%%%%%%%%%%%%%%%%%%%%%%%%
%%%%%%%%%%%%%%%%%%%%%%%%%%%%%%%%%%%%%%%%%%%%%%%%%%%%%%%%%%%%
By triangular inequality for integrals, we get
\begin{equation}\label{Eq16B}
\begin{split}
I^{NL} &\lesssim  \sum_{(k_{0}, k_{1}, k_{2}) \in \Theta(N) } \frac{N_{0}^{s}}{ N_{1}^{s}  N_{2}^{s}}
 \frac{1}{L_{0}^{b'}L_{1}^{b} L_{2}^{b}} \times \\
& \qquad \times \left|\int_{ D} \left( \int_{M} \tilde{P}_{k_{0}} \overline{\widehat{w_{0}}(\tau_{0})} P_{k_{1}}\widehat{w_{1}}(\tau_{1}) \tilde{P}_{k_{2}} \overline{\widehat{w_{2}}(\tau_{2})} dg \right) \delta(\tau_{1} - \tau_{2} - \tau_{0}) d \tau_{0} d \tau_{1} d \tau_{2}   \right| \\
 &\lesssim  \frac{N_{0}^{s}}{ N_{1}^{s}  N_{2}^{s}}
 \frac{1}{L_{0}^{b'}L_{1}^{b} L_{2}^{b}} \sum_{(k_{0}, k_{1}, k_{2}) \in \Theta(N) }
  \int_{D} \left|\int_{M} P_{k_{0}} \widehat{w_{0}}(\tau_{0}) \tilde{P}_{k_{1}} \overline{\widehat{w_{1}}(\tau_{1})} P_{k_{2}} \widehat{w_{2}}(\tau_{2}) dg  \right| d \mu,
\end{split}
\end{equation}
where  $D$ is equipped with the measure  $d\mu := \delta(\tau_{1} - \tau_{2} - \tau_{0}) d \tau_{0} d \tau_{1} d \tau_{2}$, and 
\[
\Theta(N) := \{(k_{0}, k_{1}, k_{2}) \mid  N_{j} \leq \langle \mu_{k_{j}} \rangle^{\frac{1}{2}} < 2 N_{j}, \; j = 0,1, 2 \} .
\]

By Lemma \ref{EqStar}, and noting that $\mu_{k_{0}}^{-p} \leq 2^{p} N_{0}^{-p}$ from \eqref{Eq16B}, we obtain
\begin{equation}\label{Eq18}
\begin{split}
 I^{NL} &\lesssim_{p}  \frac{N_{0}^{s}}{ N_{1}^{s}  N_{2}^{s}}
 \frac{1}{L_{0}^{b'}L_{1}^{b} L_{2}^{b}} \sum_{(k_{0}, k_{1}, k_{2}) \in \Theta(N) }
  \int_{D}  \mu_{k_{0}}^{-p}\prod_{j = 0}^{2} \|\widehat{w_{j}}(\tau_{j})\|_{L^{2}(M)} d\mu \\
 &\lesssim_{p}  \frac{N_{0}^{s}}{ N_{1}^{s}  N_{2}^{s}} C_{p} N_{0}^{-p}
 \frac{1}{L_{0}^{b'}L_{1}^{b} L_{2}^{b}} \sum_{(k_{0}, k_{1}, k_{2}) \in \Theta(N) }
  \int_{D}  \prod_{j = 0}^{2} \|\widehat{w_{j}}(\tau_{j})\|_{L^{2}(M)} d\mu.
\end{split}
\end{equation}
%%%%%%%%%%%%%%%%%%%%%%%%%%%%%%%%%%%%%%%%%%%%%%%%%%%%%%%%%%%%
%%%%%%%%%%%%%%%%%%%%%%%%%%%%%%%%%%%%%%%%%%%%%%%%%%%%%%%%%%5
Let
\begin{equation}\label{Eq19}
I_{D} := \int_{D}\prod_{j = 0}^{2} \|\widehat{w_{j}}(\tau_{j})\|_{L^{2}(M)} \delta(- \tau_{0} +\tau_{1} - \tau_{2}) d \tau_{0} d \tau_{1} d \tau_{2}.
\end{equation}

Performing the change of variables $(\tau_0, \tau_1) \mapsto (\tau_0', \tau_1' + \tau_0')$ one may rewrite \eqref{Eq19} as 
\begin{equation}\label{Eq19 A}
I_{D'} = \int_{D'} \|\widehat{w_{0}}(\tau_{0}')\|_{L^{2}(M)} \|\widehat{w_{1}}(\tau_{1}' + \tau_{0}')\|_{L^{2}(M)} \|\widehat{w_{2}}(\tau_{1}')\|_{L^{2}(M)} d \tau_0' d \tau_1'
\end{equation}
with, 
\[
D' = \Big\{ (\tau_0', \tau_1') \in \mathbb{R}^2: \langle \tau_0' + \mu_0 \rangle \sim L_0, \langle \tau_0' + \tau_1' + \frac{\mu_{k_1}}{\sigma} \rangle \sim L_1 \mbox{ and } \langle \tau_1' + \mu_{k_2} \rangle \sim L_2 \Big \}, 
\]
where we used the fact that for  $f$ measurable and $x_0 \in \mathbb{R}$, we have 
\[
\int_{\mathbb R} f(x) \delta_{x_0}(x)  dx  = f(x_0).  
\]
Using the Cauchy-Schwarz inequality with respect to $\tau_0'$ we obtain  
\begin{equation}
\begin{split}
I_{D'} &\leq \int_{\tau_1'}  \left(\int_{\tau_0'} \|\widehat{w_{0}}(\tau_{0}')\|_{L^{2}(M)} d \tau_0'\right)^{1/2} \left(\int_{\tau_0'}\|\widehat{w_{1}}(\tau_{1}' + \tau_{0}')\|_{L^{2}(M)} d \tau_0'\right)^{1/2} d \tau_1' \\
&\leq \|w_{0}\|_{L^{2}( \mathbb{R} \times M)}  \|w_{1}\|_{L^{2}( \mathbb{R} \times M)}  \int_{\tau_1'} \|\widehat{w_{2}}(\tau_{1}')\|_{L^{2}(M)} d \tau_1'    
\end{split}    
\end{equation}
Now, applying  the Cauchy-Schwarz inequality with respect to $\tau_1'$ we get
\[
I_{D'} \leq \|w_{0}\|_{L^{2}( \mathbb{R} \times M)}  \|w_{1}\|_{L^{2}( \mathbb{R} \times M)}  \|w_{2}\|_{L^{2}( \mathbb{R} \times M)} \left(\int_{\tau_1' \in D'_1 } d \tau_1'\right)^{1/2}. 
\]
where 
\[
D_1' := \Big\{ \tau_1' : \langle  \tau_1' + \tau_0' + \frac{\mu_{k_1}}{\sigma} \rangle \sim L_1 \mbox{ and } \langle \tau_1' + \mu_{k_2} \rangle \sim L_2 \Big \}.
\]
Since 
\[
|D_1'| \lesssim (\min \{L_1 , L_2\})^{1/2},
\]
one obtains 
\begin{equation}\label{Eq17B}
 I_{D} \lesssim     \|w_{0}\|_{L^{2}( \mathbb{R} \times M)}  \|w_{1}\|_{L^{2}( \mathbb{R} \times M)}  \|w_{2}\|_{L^{2}( \mathbb{R} \times M)}  L_2^{1/2}. 
\end{equation}

%%%%%%%%%%%%%%%%%%%%%%%%%%%%%%%%%%%%%%%%%%%%%%%%%%%%%%

Thus, from  \eqref{Eq18} and \eqref{Eq17B}, we have
\begin{equation}\label{Eq20b}
\begin{split}
 I^{NL} & \lesssim_{p}   \frac{N_{0}^{s}}{ N_{1}^{s}  N_{2}^{s}}  N_{0}^{-p}
 \frac{1}{L_{0}^{b'}L_{1}^{b} L_{2}^{b}} \sum_{(k_{0}, k_{1}, k_{2}) \in \Theta(N) } I_{D}  \\
 & \lesssim   \frac{N_{0}^{s}}{ N_{1}^{s}  N_{2}^{s}} N_{0}^{-p}
 \frac{1}{L_{0}^{b'}L_{1}^{b} L_{2}^{b}} L_{2}^{\frac{1}{2}} \prod_{j = 0}^{2} \|w_{j}\|_{L^{2}(\mathbb{R} \times M)} \cdot \sharp \Theta(N) \\
& \lesssim_{\sigma,p}   \frac{N_{0}^{s}}{ N_{1}^{s}  N_{2}^{s}} N_{0}^{-p}
 \frac{L_{2}^{\frac{1}{2}} }{L_{0}^{b'}(L_{1} L_{2})^{b}}   (N_{0} N_{1}N_{2})^{2} \prod_{j = 0}^{2} \|w_{j}\|_{L^{2}(\mathbb{R} \times M)}.
\end{split}
\end{equation}
In the last line of \eqref{Eq20b} we used the Weyl's Law for eigenvalues of Laplace operator  in compact manifolds to find an estimate for the number of elements of the set $\Theta (N)$, that is,
\[
 \sharp \{ \mbox{ eigenvalues of   } -\Delta \leq A \} \sim \frac{\omega_{d} Vol_{g}(M)}{(2 \pi)^{d}} A^{\frac{d}{2}}.
\]

Suppose  $N_{0} > C (N_{1}+ N_{2})$. Thus, from  \eqref{Eq20b},  we obtain
\[
 I^{NL} \lesssim_{ p}   \frac{N_{0}^{s + 6 - p}}{ N_{1}^{s}  N_{2}^{s}}
 \frac{1}{L_{0}^{b'}L_{1}^{b} L_{2}^{b- \frac{1}{2}}}    \prod_{j = 0}^{2} \|w_{j}\|_{L^{2}(\mathbb{R} \times M)}.
\]
By Lemma \ref{EqStar}, we can choose $p > 0$ such that $s + 6 - p < 0$. Thus,
\[
 I^{NL} \lesssim   \frac{N_{0}^{s + 6 - p }}{ N_{1}^{s}  N_{2}^{s}}
 \frac{1}{L_{0}^{b'}L_{1}^{b} L_{2}^{b-\frac{1}{2}}}    \prod_{j = 0}^{2} \|w_{j}\|_{L^{2}(\mathbb{R} \times M)},
\]
as desired. 
\end{prova}

Since we will be dealing with dyadic sums involving different regimes on dyadic variables, the following lemma will be useful, whose proof can be found in \cite{BGT1}, page 282.
\begin{lema} \label{DyaSum}(Discrete Schur lemma). 
For all $\gamma > 0$ and all $\theta > 0$ there is  $C > 0$ such that if   $(c_{N})$ and   $(d_{N'})$
are two sequences of non-negative numbers indexed by  dyadic integers, then
\begin{equation}\label{DyaSumA}
\sum_{N \leq \gamma N'} \left(\frac{N}{N'}\right)^{\theta} c_{N} d_{N'} \leq C \left(\sum_{N} c_{N}^{2}\right)^{\frac{1}{2}}  \left(\sum_{N'} d_{N'}^{2}\right)^{\frac{1}{2}}.
\end{equation}
\end{lema}

\subsection{Proof of bilinear estimates}
We are now ready to provide proofs of the  bilinear estimates related to
quadratic interactions stated  in Proposition \ref{Bili1}.\\

\noindent
\begin{dem} [Proof of Proposition \ref{Bili1}] 
Let us prove first \eqref{Bili1a}. By density, we can assume that    $u_ {j} \in C_{0}^{\infty} (\mathbb {R} \times M) $ ($ j = 1, 2 $). As in \cite{BGT2}, we introduce  the  functions
\[
w_{1}(t) = \frac{1}{2\pi} \sum _{k \in \mathbb{N}} \langle \mu_{k}\rangle^{\frac{s}{2}} \int_{\mathbb{R}}  \left\langle \tau + \frac{\mu_{k}}{\sigma}\right\rangle^{b} \widehat{P_{k} u_{1}}(\tau) e^{i t \tau} d\tau
\]
\[
w_{2}(t) = \frac{1}{2\pi} \sum _{k \in \mathbb{N}} \langle \mu_{k}\rangle^{\frac{s}{2}} \int_{\mathbb{R}}  \langle \tau + \mu_{k}\rangle^{b} \widehat{P_{k} u_{2}}(\tau) e^{i t \tau} d\tau.
\]

Using Proposition \ref{Eq01}, we can write
\[
u_{1}(t) = \frac{1}{2\pi} \sum _{k \in \mathbb{N}} \langle \mu_{k}\rangle^{-\frac{s}{2}} \int_{\mathbb{R}}  \left\langle \tau + \frac{\mu_{k}}{\sigma}\right\rangle^{-b} \widehat{P_{k} w_{1}}(\tau) e^{i t \tau} d\tau
\]
and 
\[
u_{2}(t) = \frac{1}{2\pi} \sum _{k \in \mathbb{N}} \langle \mu_{k}\rangle^{-\frac{s}{2}} \int_{\mathbb{R}}  \langle \tau + \mu_{k}\rangle^{-b} \widehat{P_{k} w_{2}}(\tau) e^{i t \tau} d\tau,
\]
with 
\begin{equation}\label{Eqq1}
\|u_{1}\|_{X_{1/\sigma}^{s, b}} = \|w_{1}\|_{L^{2}(\mathbb{R} \times M)} \mbox{ } \mbox{ and  } \mbox{ } 
\|u_{2}\|_{X^{s, b}} = \|w_{2}\|_{L^{2}(\mathbb{R} \times M)}.
\end{equation}

Using the duality relation between  $X^{s, -b'}$ and $ X^{-s, b'} \approx (X^{s, -b'})^{\ast} $ and \eqref{Eq1}, it follows that to prove \eqref{Bili1a}, it suffices to show that
\begin{equation}\label{Eqq0}
I := \left|\int_{\mathbb{R} \times M} \overline{u_0} u_1 \overline{u_2}  \right| \leq C   \|w_1\|_{L^{2}(\mathbb{R} \times M)} \|w_2\|_{L^{2}(\mathbb{R} \times M)}  \|w_0\|_{L^{2}(\mathbb{R} \times M)},
\end{equation}
where  $w_{0} \in L^{2}(\mathbb{R} \times M)$ is arbitrary, and 
\[
u_0(t) = \frac{1}{2\pi} \sum _{k \in \mathbb{N}} \langle \mu_{k}\rangle^{\frac{s}{2}} \int_{\mathbb{R}}  \langle \tau + \mu_{k}\rangle^{-b'} \widehat{P_k w_0}(\tau) e^{i t \tau} d\tau.
\]

Now, we move to estimate   the term $I$ in \eqref{Eqq0}. Let $N_0, N_1, N_2 $ and $ L_0, L_1, L_2 $ be the dyadic integers, that is,  $N_{j} = 2^{n_j}$ $(n_j \in \mathbb{N})$, $L_j = 2^{\ell_j}$ $(\ell_j \in \mathbb{N})$, $j = 0, 1,2$.
We define $ N: = (N_0, N_1, N_2) $, $ L: = (L_0, L_1, L_2) $, and  use this notation throughout this proof. The sum $ \Sigma_{N}$
denotes  the summation over all possible dyadic values of               $ N_0, N_1, N_2 $. Similar convention will be adopted for the sum over $ L $. 

Note that, we can decompose\footnote {See subsection  \ref{projetors} for more details on these spectral decompositions.}
the functions  $u_{j}$ is a such way that 
\[
u_j(t, \cdot) =  \sum_{N_j} \sum_{L_j } u_j^{N_j L_j}(t, \cdot),
\] 
with respective components 
\begin{equation}\label{Eq12}
u_j^{N_j L_j}(t) = \frac{1}{2\pi} \sum_{ k: N_j \leq \langle \mu_{k}\rangle^{1/2} < 2 N_j} \langle \mu_{k}\rangle^{\pm s/2}  \int_{L_j \leq \langle \tau + \tilde{\mu_k}\rangle < 2 L_j} \langle \tau + \widetilde{\mu_k}\rangle^{- \tilde{b}} \widehat{P_k w_j}(\tau) e^{i t \tau} d\tau,
\end{equation}
 where $\widetilde{\mu_k} = \frac{\mu_k}{\sigma} $  and  $\tilde{b} = b$ or  $b'$. 
 
 Inserting this decomposition in the LHS of  \eqref{Eqq0}, one obtains
\begin{equation}\label{Eq12a}
I \lesssim  \sum_{L} \sum_{N} I^{NL},
\end{equation}
where 
\begin{equation}\label{Eq12b}
I^{NL} := \left|\int_{\mathbb{R} \times M} u_{0}^{N_{0}L_{0}} \overline{u_{1}^{N_{1}L_{1}}} u_{2}^{N_{2}L_{2}}  \right|,
\end{equation}
with 
\begin{equation}\label{Eq12c}
u_{0}^{N_{0}L_{0}}(t) = \frac{1}{2\pi} \sum_{N_{0} \leq \langle \mu_{k}\rangle < 2 N_{0}} \langle \mu_{k}\rangle^{s/2}  \int_{ L_{0}\leq \langle \tau + \mu_{k}\rangle <  2 L_0} \langle \tau + \mu_{k}\rangle^{-b'}  \widehat{P_{k} w_{0}}(\tau) e^{i t \tau} d\tau,
\end{equation}
\begin{equation}\label{Eq12d}
u_{1}^{N_{1}L_{1}}(t) = \frac{1}{2\pi} \sum_{N_{1} \leq \langle \mu_{k}\rangle < 2 N_{1}} \langle \mu_{k}\rangle^{-s/2}  \int_{ L_{1}\leq \langle \tau + \frac{\mu_{k}}{\sigma}\rangle <  2 L_{1}} \langle \tau + \frac{\mu_{k}}{\sigma} \rangle^{-b}  \widehat{P_{k} w_{1}}(\tau) e^{i t \tau} d\tau,
\end{equation}
and 
\begin{equation}\label{Eq12e}
u_{2}^{N_{2}L_{2}}(t) = \frac{1}{2\pi} \sum_{N_{2} \leq \langle \mu_{k}\rangle < 2 N_{2}} \langle \mu_{k}\rangle^{-s/2}  \int_{ L_{2}\leq \langle \tau + \mu_{k}\rangle <  2 L_{2}} \langle \tau + \mu_{k}\rangle^{-b}  \widehat{P_{k} w_{2}}(\tau) e^{i t \tau} d\tau.
\end{equation}

Applying  the Fubini's theorem, we obtain from \eqref{Eq12} that
\begin{equation}\label{Eq13}
\begin{split}
 \frac{1}{2\pi}  \int_{\langle \tau + \tilde{\mu_{k}}\rangle \sim L_{j}} \langle \tau + \tilde{\mu_{k}}\rangle^{-\tilde{b}} \widehat{P_{k} w_{j}}(\tau) e^{i t \tau} d\tau  &=  \frac{1}{2\pi}  \int_{\mathbb{R}} \chi_{ \{\tau: \langle \tau + \tilde{\mu_{k}}\rangle \sim L_{j} \}} (\tau) \langle  \tau + \tilde{\mu_{k}}\rangle^{-\tilde{b}} \widehat{P_{k} w_{j}}(\tau) e^{i t \tau} d\tau  \\
  &=  \frac{1}{2\pi}  \int_{\mathbb{R}} \chi_{ \{\tau: \langle \tau + \tilde{\mu_{k}}\rangle \sim L_{j} \}} (\tau) \langle  \tau + \tilde{\mu_{k}}\rangle^{-\tilde{b}} P_{k} \widehat{ w_{j}}(\tau) e^{i t \tau} d\tau  \\
 &  = P_{k} \Big( \chi_{ \{\tau: \langle \tau + \tilde{\mu_{k}} \rangle \sim L_{j}  \}}(\tau) \langle \tau + \tilde{\mu_{k}}\rangle^{- \tilde{b}} \widehat{ w_{j}}(\tau) \Big)^{\vee}(t). \\
\end{split}
\end{equation}

Inserting  \eqref{Eq13} in  \eqref{Eq12}, and taking the  $L^{2}(\mathbb{R} \times M)$-norm, we get
\begin{equation}\label{Eq14}
\begin{split}
\|u_{j}^{N_{j}L_{j}}\|_{L^{2}(\mathbb{R} \times M)}^{2} & =  \sum_{\langle \mu_{k}\rangle \sim N_{j}^{2}} \langle \mu_{k}\rangle^{\pm s}  \int_{\langle \tau + \tilde{\mu_{k}}\rangle \sim L_{j}} \langle \tau + \tilde{u_{k}}\rangle^{- 2 \tilde{b}} \|\widehat{P_{k} w_{j}}(\tau)\|_{L^{2}(M)}^{2}  d\tau \\
 & \lesssim N_{j}^{\pm 2s} L_{j}^{- 2 \tilde{b}} (c_{j}^{N_{j}L_{j}})^{2}, 
\end{split}
\end{equation}
where, 
\begin{equation}\label{cjCoef}
(c_{j}^{N_{j}L_{j}})^{2} := \sum_{\langle \mu_{k}\rangle \sim N_{j}^{2}}   \int_{\langle \tau + \tilde{\mu_{k}}\rangle \sim L_{j}}  \|\widehat{P_{k} w_{j}}(\tau)\|_{L^{2}(M)}^{2}  d\tau.
\end{equation}

Now, summing on the dyadic variables  $L_{j}$ and $N_{j}$, 
and using  the  Plancherel's theorem  in the time variable, we obtain
\begin{equation}\label{Eq15}
\sum_{L_j} \sum_{N_j} (c_j^{N_j L_j})^{2}  \lesssim \|w_j\|_{L^{2}(\mathbb{R} \times M)}^{2}.
\end{equation}

Considering the coefficients  $c_j^{N_j L_j}$,  we need to estimate the components of $u_j$ given in \eqref{Eq12} in the norms of the $X^{s,b}$ and  $X^{s,b}_{1/\sigma}$  spaces. For this, we  consider an appropriate selection of real parameters $a_{1}, b_{1}, a_{2}, b_{2}$, to use the estimate  \eqref{LemaApB1} of the Lemma \ref{LemaApB} and then  \eqref{Eq14} to obtain
\begin{equation}\label{Eq15a}
\|u^{N_0 L_0}_{0}\|_{X^{a_2, b_2}} \lesssim  L_{0}^{b_2} N_0^{a_2}\|u^{N_0 L_0}_{0}\|_{L^{2}(\mathbb{R} \times M)} \lesssim L_0^{b_{2} - b'} N_{0}^{a_2 + s_2}c_{0}^{L_0 N_0}, 
\end{equation}
\begin{equation}\label{Eq15b}
\|u^{N_1 L_1}_{1}\|_{X_{1/\sigma}^{a_1,b_1}} \lesssim  L_{1}^{b_1 } N_{1}^{a_1} \|u^{N_{1}L_{1}}_{1}\|_{L^{2}(\mathbb{R} \times M)} \lesssim  L_{1}^{b_{1} - b} N_{1}^{a_{1} - s_1} c_{1}^{N_{1}L_{1}},
\end{equation}
\begin{equation}\label{Eq15c}
\|u^{N_2 L_2}_{2}\|_{X^{a_1,b_1}} \lesssim  L_{2}^{b_1} N_{2}^{a_1} \|u^{N_2 L_2}_{2}\|_{L^{2}(\mathbb{R} \times M)}  \lesssim L_{2}^{b_{1} - b} N_{2}^{a_{1} - s_2} c_{2}^{N_{2} L_{2}}.
\end{equation}
By means  of the Lemma \ref{LemaFund}, using \eqref{Eq12a}, it is convenient to write
\begin{equation}
I \leq I_{1} + I_{2},
\end{equation}
where in  $I_{1}$  the summation  in $N$ is restricted to $N_{0} \geq  C (N_{1} + N_{2})$ and all other possibilities are in  $I_{2}$.

To bound $I_{1}$, note that the frequency regime over $N = (N_{0}, N_{1}, N_{2})$ is such that $N_{0} \geq C (N_{1}+ N_{2})$ so that we can apply the estimate \eqref{LemaFunda} of  Lemma \ref{LemaFund} with $ p > s + 6 $. As $ b > \frac{1}{2} $ and $ b' > 0$, we can perform a sum of geometric series\footnote {By this we mean that the sums in question are of the form $\sum_{ k=0}^{\infty} \left(\frac{1}{2^{\gamma}}\right)^{k} < \infty$, since $ \gamma > 0 $.} in all dyadic variables to conclude that

\begin{equation}\label{BoundForI1}
I_{1} \lesssim \sum_{L} \sum_{N_{0} \geq  C (N_{1} + N_{2})}  \frac{N_{0}^{s + 6 - p}}{ N_{1}^{s}  N_{2}^{s}}
 \frac{1}{L_{0}^{b'}L_{1}^{b -1/2} L_{2}^{b}}    \prod_{j = 0}^{2} \|w_{j}\|_{L^{2}(\mathbb{R} \times M)} \lesssim \prod_{j = 0}^{2} \|w_{j}\|_{L^{2}(\mathbb{R} \times M)}.
\end{equation}

 Having the estimate \eqref{BoundForI1} at hand, it remains to bound 
\begin{equation}\label{IntegralI2}
 I_{2} := \sum_{L} \sum_{N_0 < C (N_1 + N_2)} I^{NL}.
\end{equation}

In this case, since the frequency regime satisfies the condition $N_0 \leq C (N_1 + N_2)$, we cannot apply  Lemma \ref{LemaFund}. To overcome this obstacle, we consider the cases $N_{1} \geq N_{2}$ and  $N_{1} \leq N_{2}$ separately and find  appropriate estimates in each case.

 By symmetry, we just consider the case where $N_1 \geq N_2$. In this case, $N_0 \leq 2C N_1$. Once this frequency regime has been set, our next step will be to bound the term $I^{NL}$ in \eqref{IntegralI2} in two different ways. The first way is to bound $ I^{NL}$ using Cauchy-Schwarz inequality which behaves  better with respect to the localization on $ N $. The other way is to bound  $ I^{NL}$ using 
 H\"older's inequality which behaves  better with respect to the localization on $ L $. Interpolating these bounds, we obtain the required estimate.  In the sequel, we describe this process in detail.
 
An use of the Cauchy-Schwarz inequality with respect to $N$ in \eqref{Eq12b}, yields
\begin{equation}\label{Eq21}
I^{NL} \lesssim \|u_{0}^{L_0 N_0}\|_{L^{2}(\mathbb{R} \times M)}  \|\overline{u_{1}^{L_{1} N_{1}}} u_{2}^{L_{2} N_{2}}\|_{L^{2}(\mathbb{R} \times M)}.
\end{equation}
Applying \eqref {Eq14} directly, it follows that
\begin{equation}\label{Eq21b}
\|u_0^{L_0 N_0}\|_{L^{2}(\mathbb{R} \times M)} \lesssim L_0^{- b'} N_0^{s}c_0^{L_0 N_0}. 
\end{equation}

Now, using the bilinear estimate \eqref{Eq00b} with  $s' > s_0(M)$, we have
\begin{equation}\label{Eq21a}
\|\overline{u_1^{L_1 N_1}} u_2^{L_2 N_2}\|_{L^{2}(\mathbb{R} \times M)} \lesssim \min(N_1, N_2)^{s'} \|u_1^{L_1 N_1}\|_{X_{1/\sigma}^{0, \frac{1}{2} + \epsilon_0}(\mathbb{R} \times M)} \|u_2^{L_2 N_2}\|_{X^{0, \frac{1}{2} + \epsilon_0 }(\mathbb{R} \times M)}.
\end{equation}
 Considering  \eqref{Eq15b} and  \eqref{Eq15c} with the values  $a_{i} = 0$, $b_{i} = \frac{1}{2} + \epsilon_{0} \mbox{ }(i =  1, 2)$ and  $s_1 = s$,  (where $\epsilon_{0} > 0 $ and  $b > 1/2$ will be chosen later in a suitable way), we get
\begin{equation}\label{Eq21c}
\|u^{N_1 L_1}_{1}\|_{X_{1/\sigma}^{0, \frac{1}{2} + \epsilon_0 }}  \lesssim_{\sigma}  L_{1}^{ \frac{1}{2} + \epsilon_{0} - b} N_1^{- s} c_1^{L_1 N_1}, 
\end{equation}
and
\begin{equation}\label{Eq21d}
\|u^{N_2 L_2}_{2}\|_{X^{0,  \frac{1}{2} + \epsilon_0 }} \lesssim L_{2}^{\frac{1}{2} + \epsilon_{0} - b} N_{2}^{ - s_2} c_{2}^{L_{2} N_{2}}.
\end{equation}

We note here that the value of  $s_{2}$ in \eqref{Eq21d} will be chosen so that $s_{2} = 1$ if $s \geq 1$ and $s_{2} = s$ if $0 < s < 1$. 
Replacing the estimates  \eqref{Eq21b} and  \eqref{Eq21a}
in \eqref{Eq21}, we obtain
\begin{equation}\label{INLA}
I^{NL} \lesssim  N_{2}^{s'}  \frac{( L_1 L_2)^{\frac{1}{2} + \epsilon_0}}{L_0^{b'} L_1^{b} L_2^{b}} \frac{N_0^{s}}{N_1^{s} N_2^{s_2}} \prod_{j = 0}^{2} c_j^{L_j N_j}.
\end{equation}

Now, we make a better estimate with respect to localization over $L$. In fact, by applying
H\"older's inequality in \eqref{Eq12b}, we get

\begin{equation}\label{Eq23}
I^{NL} \lesssim \|u_{0}^{L_{0} N_{0}}\|_{L^{3}(\mathbb{R}; L^{2}(M))} \|u_{1}^{L_{1} N_{1}}\|_{L^{3}(\mathbb{R}; L^{2}(M))} \|u_{2}^{L_{2} N_{2}}\|_{L^{3}(\mathbb{R}; L^{\infty}(M))}.
\end{equation}
By using the Sobolev embedding  $H^{\frac{d}{2}+ \frac{1}{10}}(M) \hookrightarrow L^{\infty}(M)$ in \eqref{Eq23}, we obtain 
\begin{equation}\label{Eq23a}
 \|u_{2}^{L_{2} N_{2}}\|_{L^{\infty}(M)} \lesssim  \|u_{2}^{L_{2} N_{2}}\|_{H^{\frac{d}{2}+ \frac{1}{10}}(M)} \lesssim N_2^{\frac{d}{2}+ \frac{1}{10}}  \|u_{2}^{L_{2} N_{2}}\|_{L^{2}(M)}.
\end{equation}
Thus, from \eqref{Eq23} we have
\begin{equation}\label{Eq24}
I^{NL} \lesssim N_{2}^{\frac{d}{2}+ \frac{1}{10}}\|u_{0}^{L_{0} N_{0}}\|_{L^{3}(\mathbb{R}; L^{2}(M))} \|u_{1}^{L_{1} N_{1}}\|_{L^{3}(\mathbb{R}; L^{2}(M))} \|u_{2}^{L_{2} N_{2}}\|_{L^{3}(\mathbb{R}; L^{2}(M))}.
\end{equation}
Using part $(ii)$ of  Lemma \ref{BasicXsb}, we obtain $\|f\|_{L^{3}(\mathbb{R}; L^{2}(M))} \lesssim \|f\|_{Y}$ where $Y = X_{1/\sigma}^{0, \frac{1}{6}}(\mathbb{R} \times M)$ or $X^{0, \frac{1}{6}}(\mathbb{R} \times M)$. Thus, using the identities  involving the coefficients $c_{j}^{L_{j} N_{j}}$ given in  \eqref{Eq15a}, \eqref{Eq15b} and  \eqref{Eq15c}, it follows that
\begin{equation}\label{INLB}
\begin{split}
I^{NL} & \lesssim N_2^{\frac{d}{2}+ \frac{1}{10}}\|u_{0}^{L_0 N_0}\|_{X^{0, \frac{1}{6}}} \|u_{1}^{L_1 N_1}\|_{X_{1/\sigma}^{0, \frac{1}{6}}} \|u_{2}^{L_2 N_2}\|_{X^{0, \frac{1}{6}}}\\
 & \lesssim N_2^{\frac{d}{2}+ \frac{1}{10}} L_0^{\frac{1}{6} - b'} N_0 ^{s} L_1^{\frac{1}{6} - b} N_1 ^{-s} L_2^{\frac{1}{6} - b} N_2 ^{-s_2} \prod_{j = 0}^{2} c_j^{L_j  N_j} \\
  & = N_2 ^{\frac{d}{2}+ \frac{1}{10}}  \frac{(L_0 L_1 L_2)^{\frac{1}{6} }}{L_0^{b'} L_1^{b} L_2^{b}} \frac{N_0^{s}}{N_1^{s} N_2^{s_2}} \prod_{j = 0}^{2} c_j^{L_j N_j}. 
\end{split}
\end{equation}

 An interpolation between  \eqref{INLA} and \eqref{INLB}  implies that for every  $\theta \in (0, 1)$
\begin{equation}\label{INL0}
I^{NL} \lesssim
   \Big(\frac{N_0}{N_1} \Big)^{s} N_2^{(\frac{d}{2}+ \frac{1}{10} - s_2) \theta + (1 - \theta) (s' - s_2)} \frac{(L_0 L_1 L_2)^{\frac{\theta}{6} + (\frac{1}{2} + \epsilon_0) \theta }}{L_0^{b'} L_1^{b} L_2^{b}} \prod_{j = 0}^{2} c_{j}^{L_j N_j}.
\end{equation}

Now, we consider a fixed value of the parameter $s$ such that 
  $s' < s$,  where  $0 < s_0(M) < s' < 1$ and make an analysis dividing in  two different cases. 

\noindent
{\bf Case 1. $s \in (0, 1)$.} In this case $s_2 = s$ and  \eqref{INL0} reduces to
\begin{equation}\label{INL0A}
I^{NL} \lesssim
   \Big(\frac{N_{0}}{N_{1}} \Big)^{s} N_{2}^{(\frac{d}{2}+ \frac{1}{10} - s) \theta + (1 - \theta) (s' - s)} \frac{(L_{0} L_{1} L_{2})^{\frac{\theta}{6} + (\frac{1}{2} + \epsilon_0) \theta }}{L_{0}^{b'} L_{1}^{b} L_{2}^{b}} \prod_{j = 0}^{2} c_j^{L_j N_j}.
\end{equation}

In order to make the exponent of $N_2$ negative (necessary for convergence), we  choose $ \theta $ so that
$0 < \theta < \frac{s - s'}{\frac{d}{2}+ \frac{1}{10} - s'} < 1$. Now, we choose $\epsilon_0$
such that $0 < \epsilon_0 < \frac{\theta}{9 - 3 \theta}$. Hence, if we define $b' := \frac{1}{2} - 2 \epsilon_0$, this ensures that
 $ b' > \frac{\theta}{6} + (\frac{1}{2} + \epsilon_0) \theta$. Moreover, one has $b' > 0 \Longleftrightarrow \epsilon_0 < \frac{1}{4}$, which is a consequence of $\frac{\theta}{9 - 3 \theta} \leq \frac{1}{4} \Longleftrightarrow \theta \leq \frac{9}{7} = 1 + \frac{2}{7} $. Finally, choose  $b : = \frac{1}{2} + \frac{3}{2} \epsilon_0$ and notice that
 with this choice of parameters, the basic conditions for $(b, b')$ are verified, i.e., 
 $0< b' < \frac{1}{2} < b $ and $b + b' < 1 $. 

\noindent{\bf Case 2. $s \in [1, + \infty)$.} In this case $s_2 = 1$ and  \eqref{INL0}  is reduced to 
\begin{equation}\label{INL0B}
I^{NL} \lesssim
   \Big(\frac{N_{0}}{N_{1}} \Big)^{s} N_{2}^{ (\frac{d}{2}+ \frac{1}{10} -1)\theta  + (1 - \theta) (s' -1)} \frac{(L_{0} L_{1} L_{2})^{\frac{\theta}{6} + (\frac{1}{2} + \epsilon_0) \theta }}{L_{0}^{b'} L_{1}^{b} L_{2}^{b}} \prod_{j = 0}^{2} c_{j}^{L_{j}N_{j}}.
\end{equation}

Observe that
\[
 \Big(\frac{d}{2}+ \frac{1}{10} -1 \Big)\theta  + (1 - \theta) (s' -1) = \theta \Big(\frac{d}{2}+ \frac{1}{10} - s'\Big) + (s' -1).
\]
As $0 <  s' < 1$,  the choice of  $\theta$ to obtain the exponent of
 $N_{2}$ negative is 
\[
0 < \theta < \frac{1 - s'}{\frac{d}{2}+ \frac{1}{10} - s'} < 1. 
\]

Therefore, at this point, we can choose the same values of
 $\epsilon_0, b'$ and $b$ of the {\bf Case 1}.

Hence,  after the suitable choice of the parameters in 
 \eqref{INL0A}  and  \eqref{INL0B}, as $b' < b$, we have that in any case, there are $\gamma_{1},\cdots, \gamma_{4} > 0$ such that
\begin{equation}\label{INL1}
I^{NL} \lesssim
   \Big(\frac{N_{0}}{N_{1}} \Big)^{s} \frac{1 }{N_{2}^{\gamma_{1}}} \frac{1}{L_{0}^{\gamma_{2}} L_{1}^{\gamma_{3}} L_{2}^{\gamma_{4}}} \prod_{j = 0}^{2} c_{j}^{L_{j}N_{j}}.
\end{equation}

Next, using  \eqref{INL1}, it follows from  \eqref{IntegralI2} that
\begin{equation}\label{INL2}
I_{2} \lesssim
  \sum_{L} \sum_{N: N_{0} \leq 2C N_{1} }  \Big(\frac{N_{0}}{N_{1}} \Big)^{s} \frac{1 }{N_{2}^{\gamma_{1}}} \frac{1}{L_{0}^{\gamma_{2}} L_{1}^{\gamma_{3}} L_{2}^{\gamma_{4}}} \prod_{j = 0}^{2} c_{j}^{L_{j}N_{j}}.
\end{equation}
The summations involving $L_0, L_1, L_2, N_2$ in \eqref{INL1} 
can be performed via convergence of geometric series. Now, we indicate the coefficients associated with the dyadic variables
$N_{1}$, $N_{0}$ and $N_{2}$ by
\begin{equation}\label{INL3}
(\alpha_{N_{1}})^{2} = \sum_{L_{1}} (c_{1}^{L_{1}N_{1}})^{2}, \qquad  (\beta_{N_0})^{2} = \sum_{L_0} (c_0^{L_0 N_0})^{2} \qquad \mbox{ and } \qquad (\gamma_{N_2})^{2} = \sum_{L_2} (c_{2}^{L_{2} N_{2}})^{2}.
\end{equation}
Applying the Cauchy-Schwarz inequality on the variables 
$L_{j}$ ($ j = 0, 1, 2$) and using the relations given in \eqref{INL3}, we obtain
\begin{equation}\label{INL6}
\begin{split}
   I_{2} & \lesssim \sum_{(N_0, N_1): N_0 \leq 2 C N_1 }  \Big(\frac{N_0}{N_1} \Big)^{s}  \sum_{N_{2}} \frac{1 }{N_{2}^{\gamma_{1}}} \left(\sum_{L_0, L_1, L_2}  \frac{1}{L_{0}^{\gamma_{2}} L_{1}^{\gamma_{3}} L_{2}^{\gamma_{4}}} \prod_{j = 0}^{2} c_{j}^{L_{j} N_{j}} \right) \\
    & \lesssim  \sum_{(N_0, N_1): N_{0} \leq 2 C N_1 }  \Big(\frac{N_0}{N_1} \Big)^{s}  \alpha_{N_1} \beta_{N_0} \sum_{N_2} \frac{1 }{N_2^{\gamma_{1}}}  \gamma_{N_2}.  \\
\end{split}
\end{equation}

 Applying the Cauchy-Schwarz inequality on the summation involving      $N_{2}$ and using the relation \eqref{Eq15}, it follows from   \eqref{INL6} that
\begin{equation}\label{INL7}
\begin{split}    
  I_{2} & \lesssim \sum_{(N_{0}, N_{1}): N_{0} \leq 2 C N_{1}}  \Big(\frac{N_{0}}{N_{1}} \Big)^{s}  \alpha_{N_{1}} \beta_{N_{0}} \left(\sum_{N_{2}} \gamma_{N_{2}}^{2}\right)^{1/2} \\
    & \lesssim  \|w_{2}\|_{L^{2}({\mathbb{R} \times M})}\sum_{(N_{0}, N_{1}): N_{0} \leq 2 C N_{1}}  \Big(\frac{N_{0}}{N_{1}} \Big)^{s}  \alpha_{N_{1}} \beta_{N_{0}}.
\end{split}
\end{equation}

Now, from  Lemma  \ref{DyaSum},  it follows that
\begin{equation}\label{INL5}
\sum_{N_{0} \leq 2 C N_{1}} \Big(\frac{N_{0}}{N_{1}} \Big)^{s} \alpha_{N_{1}} \beta_{N_{0}} \leq \|w_{0}\|_{L^{2}(\mathbb{R} \times M)} \|w_{1}\|_{L^{2}(\mathbb{R} \times M)}.
\end{equation}

Thus, replacing  \eqref{INL5} in \eqref{INL7}, we get
\begin{equation}\label{I1-1}
  I_{2} \lesssim \|w_{0}\|_{L^{2}({\mathbb{R} \times M})} \|w_{1}\|_{L^{2}({\mathbb{R} \times M})}  \|w_{2}\|_{L^{2} ({\mathbb{R} \times M})}.
\end{equation}

Combining the estimates  \eqref{BoundForI1} and \eqref{I1-1} we obtain the required bilinear estimates \eqref{Bili1a} and \eqref{Bili1b}.

 The proof of the bilinear estimate \eqref{Bili2a} is similar to the proof of the estimate \eqref{Bili1a}
modulo some modifications. Let us highlight the main modifications here.

By using the duality relation between  $X_{1/\sigma}^{s, -b_1'}$ and $ X_{1/\sigma}^{-s, b_1'} \approx (X_{1/\sigma}^{s, -b_1'})^{\ast} $  we see that to prove \eqref{Bili2a}, it is necessary to prove that
\begin{equation*}
\tilde{I}:= \left|\int_{\mathbb{R} \times M} \overline{u_{0}} u_{1} u_{2}  \right| \leq C   \|u_{1}\|_{X^{s,b_1}(\mathbb{R} \times M)} \|u_{2}\|_{X^{s,b_1}(\mathbb{R} \times M)}  \|u_{0}\|_{X_{1/\sigma}^{-s, b_1'}(\mathbb{R} \times M)},
\end{equation*}
where $u_{0} \in X_{1/\sigma}^{-s, b_1'}(\mathbb{R} \times M)$ is  arbitrary. In this case, as in  \eqref{Eq12b}, it is necessary 
to work with 
\begin{equation*}
\tilde{I}^{NL} := \left|\int_{\mathbb{R} \times M} \overline{u_{0}^{N_{0}L_{0}}} u_{1}^{N_{1}L_{1}} u_{2}^{N_{2}L_{2}}  \right|.
\end{equation*}

Notice that the analogue of \eqref{Eq12} is given by
\begin{equation*}
\tilde{I}^{NL} \lesssim \|u_{0}^{L_{0} N_{0}}\|_{L^{2}(\mathbb{R} \times M)}  \|u_{1}^{L_{1} N_{1}} u_{2}^{L_{2} N_{2}}\|_{L^{2}(\mathbb{R} \times M)},
\end{equation*}
where we do not have the presence of conjugates.
Using the bilinear estimate \eqref{Eq00b} with $\sigma = 1$, we obtain 
\begin{equation*}
\|u_{1}^{L_{1} N_{1}}u_{2}^{L_{2} N_{2}}\|_{L^{2}(\mathbb{R} \times M)} \lesssim \min(N_{1}, N_{2})^{s} \|u_{1}^{L_{1} N_{1}}\|_{X^{0, \frac{1}{2} + \epsilon_0'}(\mathbb{R} \times M)} \|u_{2}^{L_{2} N_{2}}\|_{X^{0, \frac{1}{2} + \epsilon_0' }(\mathbb{R} \times M)}.
\end{equation*}

In the same way, we can obtain the estimate \eqref{Eq24}. The other estimates can be obtained analogously.
\end{dem}

\section{Local  and Global Theory}\label{SectionEC4}
In this section we provide proofs of  the main local and global well-posedness results of this work. 

\subsection{Local Theory}\label{SSE2C4}

We begin by recording  some basic estimates that will be important in order to prove the  well posedness result. Recall the class of the spaces $X_{\delta, \gamma}^{s, b}$ introduced in the Definition \ref{D2a} and the restriction spaces given in Definition \ref{D5}. 

\begin{prop}\label{D6} (Linear estimates in the spaces $X^{s,b}_{\delta}$). Let $M$ be a compact Riemannian manifold, 
 $b,s >0$, $v_0 \in H^{s}(M)$ and $\psi \in C_0^{\infty}(\mathbb{R})$  be such that $\psi = 1$ in $[-1, 1]$. Let
$\psi_{T}:= \psi(\frac{\cdot}{T})$, then $\psi_{T} = 1$ in $[-T, T]$. Under these assumptions, we have 
\begin{equation}\label{D7}
 \Big\| \psi_{T}(t) e^{i t (\delta \Delta - \gamma )} v_0 \Big\|_{X_{\delta }^{s,b}(\mathbb{R}\times M)} \lesssim_{\psi_{T}} \|v_{0}\|_{H^{s}(M)}.
\end{equation}
 Let $0 < b' < \frac{1}{2}$ and  $0 < b < 1 - b'$. Then, for all $F \in X_{\delta}^{s, -b'}(M)$,
\begin{equation}\label{D8}
\left\| \psi_{T}(t) \int_0^{t} e^{i (t - t')(\delta \Delta - \gamma)}F(t') dt'\right\|_{X_{ \delta }^{s,b}(\mathbb{R}\times M)} \lesssim  T^{1-b-b'}  \|F\|_{X_{ \delta }^{s,-b'}(\mathbb{R}\times M)},
\end{equation}
when $0 < T \leq 1$.
\end{prop}
\noindent

Now, we are in position to supply a proof to the result stated in Theorem \ref{BCL00} that  provides the local well-posedness for the IVP \eqref{SHGSA} posed on $d$-dimensional Zoll manifolds.\\

\noindent
\begin{dem}[Proof of Theorem \ref{BCL00}] Let $M$ be a Zoll manifold of dimension $d \geq 2$, $(v_0,u_0)\in \bH^s(M)$ with $s > s_0(d)$ and $\sigma = \frac{\beta}{\theta}$ with $\beta, \theta \in \{ n^{2}: n \in \N \}$. Let  $R > 0$ and  $T > 0$  be constants to  be suitably  chosen later. Consider the spaces  $\mathcal{X}^{s,b}_{T} := X^{s,b}_{T} \times (X_{1/\sigma}^{s, b})_{T}$ endowed with norm
\[
\|(v,u)\|_{\mathcal{X}_{T}^{s,b}} := \|v\|_{X^{s,b}_{T}} + \|u\|_{(X_{1/\sigma}^{s,b})_{T}}
\] 
 and  a closed ball of radius  $R >0$ in $\mathcal{X}^{s,b}_{T}$ given by
\[
B^{R}_{T} := \{ (v,u) \in \mathcal{X}^{s,b}_{T}; \|(v,u)\|_{
\mathcal{X}^{s,b}_{T}} \leq R \}. 
\]
   For $v, u \in B_{T}^{R}$ define the operators  
\begin{equation}\label{BCL0}
\begin{cases}
\displaystyle \Phi_{1}(u,v)(t) := \psi_{1}V(t)v_0 - i \epsilon_1 \psi_{T} \int_{0}^{t} V(t - t') \overline{v}(t') u(t') dt',\\
\displaystyle \Phi_{2}(u,v)(t) := \psi_{1}U_{\sigma}(t)u_0 - i \epsilon_2 \frac{\sigma}{2} \psi_{T}\int_{0}^{t}  U_{\sigma}(t - t')v(t')^{2} dt'.
\end{cases}
\end{equation}

 For appropriate choices of the constants $R>0$ and $T>0$, we will show that the application  $(\Phi_1, \Phi_2)$ is a contraction. 
 
 Applying  the linear estimates   \eqref{D7} and \eqref{D8},  the nonlinear estimates  \eqref{Bili1a}, \eqref{Bili1b} and \eqref{Bili2a}, \eqref{Bili2b}   in \eqref{BCL0}, using  the definition of the
 $\mathcal{X}^{s,b}_{T}$ norm and noting that $u,v \in B_T^{R}$,  one obtains that

\begin{equation}\label{BCL3}
\begin{cases}
\|\Phi_{1}(u,v)\|_{X_{T}^{s,b}} \leq c_{0} \|v_{0}\|_{H^{s}(M)} + c_{1} T^{1 - b - b'} R^{2},\\
\|\Phi_{2}(u,v)\|_{(X_{1/\sigma}^{s, b})_{T}} \leq  c_{0} \|u_{0}\|_{H^{s}(M)} + c_{1} T^{1 - b - b'} R^{2} . \\
\end{cases}
\end{equation}

Considering  $\Phi(u,v):=(\Phi_{1}(u,v), \Phi_{2}(u,v))$, we see from  \eqref{BCL3} that
\begin{equation}\label{BCL3A}
\|\Phi(u,v)\|_{ \mathcal{X}^{s,b}_{T}} \leq c_{0} \|(v_0,u_0)\|_{\bH^{s}(M)} + c_{1} T^{1- b - b'} R^{2}.
\end{equation}
Choosing $R :=  2 c_{0} \|(v_0, u_0)\|_{\bH^{s}(M)}$, we obtain
\begin{equation}\label{contra-1}
\|\Phi(u,v)\|_{\mathcal{X}^{s,b}_{T}} \leq \frac{R}{2} + c_{1} T^{1- b - b'} {R}^{2}.
\end{equation}
Thus, if we choose $ 0 < T < \Big(\frac{1}{2 c_{1} R}\Big)^{\frac{1}{1 - b - b'}}$, it is easy to see from \eqref{contra-1} that  $\Phi$ maps $ B^{R}_{T}$ into itself.

Now, writing
$u \overline{v} - \tilde{u} \overline{\tilde{v}} = ( u - \tilde{u}) \overline{v} + (\overline{v} - \overline{\tilde{v}}) \tilde{u}$, with the similar calculations as above,
we easily get

\[
\|\Phi(u,v) -\Phi(\tilde{u}, \tilde{v}) \|_{\mathcal{X}_{T}^{s,b}}  \leq C T^{1 - b - b'} R \|(u, v) - (\tilde{u}, \tilde{v})\|_{\mathcal{X}_{T}^{s,b}}.
\]

Thus, choosing  $T =T(\|(v_0, u_0)\|_{\bH^s}) > 0$, in such a way that 
\begin{equation}\label{BCL5a}
C T^{1-b-b'} R \leq \frac{1}{2}   \Longleftrightarrow T \leq \left(\frac{1}{2 C R}\right)^{\frac{1}{1-b-b'}},
\end{equation}
we get
$\Phi : B_{T}^{R} \longrightarrow  B_{T}^{R}$
is a contraction. Hence using  the Banach fixed point theorem, we conclude that, there is a unique
 $(v,u) \in B_{T}^{R}$ which solves the integral system \eqref{SHGSB} for $t \in [0, T]$, with $T \leq 1$.

As $b > \frac{1}{2}$, we have the embedding  $\mathcal{X}_{T}^{s, b} \hookrightarrow C([0,T], \bH^{s})$. Hence, for 
$T \in (0,1)$ and  $t \in (0, T)$,
\[
(v,u) \in C([0, T]; \bH^{s}(M)).
\]

To end the proof of the theorem, we will show that 
\[
\Phi: (\tilde{v_0}, \tilde{u_0}) \in B\left((v_0, u_0), r\right) \subset \bH^{s} \mapsto (v,u) \in \mathcal{X}_{T}^{s, b}
\]
 is  Lipschitz continuous  for some  $r > 0$. 
 
 Given  $T' < T$ take $r> 0$ such that $$ 0 < r^{\gamma(b, b')} < \frac{1}{2 C} \Big(\frac{1}{T'} - \frac{1}{T}\Big)$$ where $\gamma(b, b') :=  \frac{1}{1-b-b'} >0$
 and  $C >0$ is given in \eqref{BCL5a}. Given $(\tilde{v_0}, \tilde{u_0})\in \bH^s(M)$, consider the respective solution     $(\tilde{v}, \tilde{u})\in \bH^s(M)$, which exists on $[0, \tilde{T}]$, where  $\tilde{T} = \tilde{T}(\|(\tilde{v_0}, \tilde{u_0}) \|) > 0$,  satisfies $ T' < \tilde{T}$. Then, both the solutions 
are well defined in  $(0, T')$. 
 
Thus, similarly to what was done to obtain the estimate
\eqref{BCL3A}, we can get 
\[
\|(v, u)-(\tilde{v}, \tilde{u})\|_{\mathcal{X}_{T'}^{s, b}} \leq C \|(v_0, u_0) - (\tilde{v_0}, \tilde{u_0})\|_{\bH^{s}} + C (T')^{1-b-b'} R^{2} \|(v, u)-(\tilde{v}, \tilde{u})\|_{\mathcal{X}_{T'}^{s, b}}.
\]
As $ C (T')^{1-b-b'} R^{2} < 1$, 
\[
\|(v, u)-(\tilde{v}, \tilde{u})\|_{\mathcal{X}_{T'}^{s, b}} \leq \tilde{C}\|(v_0, u_0) - (\tilde{v_0}, \tilde{u_0})\|_{\bH^{s}}.
\]

Now,  using the embedding $\mathcal{X}_{T'}^{s, b} \hookrightarrow C([0, T'], \bH^{s})$, we get
\[
\sup_{0 \leq t \leq T'}\|(v, u)-(\tilde{v}, \tilde{u})\|_{\bH^{s}} \leq \tilde{C}\|(v_0, u_0) - (\tilde{v_0}, \tilde{u_0})\|_{\bH^{s}},
\]
which proves the assertion.
\end{dem}

\subsection{Global Theory}\label{SSE3C4}

In this subsection, we will use the Gagliardo-Nirenberg inequality on compact Riemannian manifolds to extend the  local solution for the IVP \eqref{SHGSA} obtained in Theorem \ref{BCL00} to the  global  one in dimensions 2 and 3.  We start with the following result.
\begin{prop}\label{GNM}
Let  $M$ be a compact Riemannian manifold of dimension $d \geq2$ and $1 < p \leq 2$.  If
$1 \leq q < r < p^{\ast} = \frac{dp}{d-p}$ and $\theta = \theta(p, q, r) = \frac{d p (r - q)}{r ( q (p - d) + d p)} \in (0, 1]$, then
\begin{equation}\label{EQ29a1}
\|u\|_{L^{r}(M)}^{\frac{p}{\theta}} \leq  \left( A_{opt} \| \nabla_{g} u\|_{L^{p}(M)}^{p} + B \|u\|_{L^{p}(M)}^{p} \right) \|u\|_{L^{q}(M)}^{\frac{p (1 - \theta)}{\theta}}
\end{equation}
where
\begin{equation}\label{EQ29a2}
A_{opt} = \inf \Big\{ A \in \mathbb{R}: \mbox{ there exists } B \in \mathbb{R} \mbox{ such that } (\ref{EQ29a1}) \mbox{ is valid}  \Big\}.
\end{equation}
\end{prop}
\noindent
\begin{dem}
See \cite{CM2008}, p. 854.
\end{dem}
\begin{obs}
The explicit value of  $A_{opt}$  is known if $q > p$ and $r = \frac{p (q-1)}{p-1}$, see
\cite{CM2008} p. 853. 
\end{obs}

Now, we prove  the global well-posedness result state in Theorem \ref{BCG00}.

\noindent
\begin{dem}[Proof of Theorem \ref{BCG00}]  First, consider $s=1$. 
Let $d = 2, 3$. If we choose  $p = 2 = q$, the 
 Gagliardo-Nirenberg inequality \eqref{EQ29a1} yields
\begin{equation}\label{MP1}
\|f\|_{L^r}  \leq A_{opt}^{\frac{\theta_r}{2}} \|\nabla f\|_{L^2}^{\theta_r}  \|f\|_{L^2}^{1 -\theta_r}  + B^{\frac{\theta_r}{2}}  \|f\|_{L^2},
\end{equation}
where $\theta_r = \frac{d}{2} - \frac{d}{r}$,  $2 < r < \infty $ if $d = 2$ and $2 < r <  \frac{2d}{d-2}$ if $d \geq 3$. Now, using  \eqref{ENERGY2}, we obtain
$$ \mathcal{E}(t) = \|\nabla v(t)\|_{L^2}^{2} + \|\nabla u(t)\|_{L^2}^{2} +
\| v(t)\|_{L^2}^{2}  + \alpha \|\nabla u(t)\|_{L^2}^{2} + \epsilon_1 \int_M \Ree(v^2(t) \overline{u(t)}) dg =  \mathcal{E}(0). $$
Hence
$$ \|(v, u)\|_{H^1}^{2} \leq |1 - \alpha| \| u\|_{L^2}^{2} + \int_M |v|^2 |u| dg + |\mathcal{E}(0)|. $$

Using the  Cauchy-Schwarz inequality, we obtain
\begin{equation}\label{MP2}
I(u,v):=\int_M |v|^2 |u| dg  \leq \| v(t)\|_{L^4}^{2}  \| u(t)\|_{L^2}.
\end{equation}
From \eqref{MASS1}, we conclude that
\[
\| u(t)\|_{L^2}^{2} \leq \frac{\mathcal{M}(0)}{2 \sigma} 
\qquad \mbox{ and } \qquad
\|v(t)\|_{L^2}^{2} \leq \mathcal{M}(0). 
\]

Using \eqref{MP1} with $r = 4$, we get
$$ \|v\|_{L^4}^{2}  \leq 2 A_{opt}^{\frac{d}{4}} \|\nabla v\|_{L^2}^{\frac{d}{2}}  \|v\|_{L^2}^{2 - \frac{d}{2}}  + 2 B^{\frac{d}{4}}  \|v\|_{L^2}^2.$$
If  $d = 2, 3$, we can  bound  \eqref{MP2} by
\begin{equation}\label{MP3}
 I(u,v) \leq  2 A_{opt}^{\frac{d}{4}} \|\nabla v\|_{L^2}^{\frac{d}{2}}  \|v\|_{L^2}^{2 - \frac{d}{2}} \| u\|_{L^2}  + 2 B^{\frac{d}{4}}  \|v\|_{L^2}^2 \| u\|_{L^2}.
\end{equation}

 From \eqref{MP3}, we obtain
\[
 \|(v, u)\|_{\bH^1}^{2} \leq |1 - \alpha| \| u\|_{L^2}^{2} + 2  \|v\|_{L^2}^{2 - \frac{d}{2}} \| u\|_{L^2} A_{opt}^{\frac{d}{4}}  \|\nabla v\|_{L^2}^{\frac{d}{2}}  + 2 B^{\frac{d}{4}}  \|v\|_{L^2}^2 \| u\|_{L^2}  +  |\mathcal{E}(0)|. 
 \]
Observe that 
\[
 2 \|v\|_{L^2}^{2 - \frac{d}{2}} \|u\|_{L^2}  = 2 (\|v\|_{L^2}^2 )^{1 - \frac{d}{4}} \|u\|_{L^2}  \leq  2 \mathcal{M}(0)^{1 - \frac{d}{4}} \sqrt{\frac{\mathcal{M}(0)}{2 \sigma} } = \sqrt{\frac{2}{\sigma}} \mathcal{M}(0)^{\frac{6 - d}{4}}.
\]
Denote
\[
C_{\sigma, A_{opt}, d } :=  A_{opt}^{\frac{d}{4}} \sqrt{\frac{2}{\sigma}} \mathcal{M}(0)^{ \frac{6 - d}{4} }.
\]
Thus, 
$$ \|(v, u)\|_{\bH^1}^{2} \leq |1 - \alpha| \| u\|_{L^2}^{2} +  C_{\sigma, A_{opt}, d }  \|(v, u)\|_{\bH^1}^{\frac{d}{2}}    + 2 B^{\frac{d}{4}}  \|v\|_{L^2}^2 \| u\|_{L^2}  +  |\mathcal{E}(0)|. $$
By Young's inequality, one obtains
$$ \left( 1 - \frac{d}{4}\right)\|(v, u)\|_{\bH^1}^{2} \leq |1 - \alpha| \| u\|_{L^2}^{2} +  (C_{\sigma, A_{opt}, d })^{ \frac{4}{4 - d}} \left(\frac{4 - d}{4}\right)  + 2 B^{\frac{d}{4}}  \|v\|_{L^2}^2 \| u\|_{L^2}  +  |\mathcal{E}(0)|. $$
That is, 
\begin{equation}\label{MP5}
\begin{split}
\|(v(t), u(t))\|_{\bH^1}^{2} &\leq  |1 - \alpha| \left(\frac{4 - d}{4}\right)^{-1}   \frac{\mathcal{M}(0)}{2 \sigma} +  (C_{\sigma, A_{opt}, d })^{ \frac{4}{4-d}}\\
& \qquad + \left(\frac{4 - d}{4}\right)^{-1}  \left( B^{\frac{d}{4}}  \sqrt{\frac{2}{\sigma}} \mathcal{M}(0)^{\frac{3}{2}}  +  |\mathcal{E}(0)|\right).
\end{split}
\end{equation}

Now we use an {\em a priori}  estimate  given by \eqref{MP5} combined with a standard argument to prove that  the solution $(v, u) \in C([0, T^{\ast}); \bH^{1}(M))$  is in fact  global in $\bH^{1}$ for  dimensions  $2$ and  $3$, that is, we have 
$T^{\ast} =  + \infty$. 

Now, we consider $s > 1$. From \eqref{Bili1b} and  \eqref{Bili2b}, we obtain
\begin{equation}\label{MP6}
\begin{cases}
\|\Phi_{1}(u,v)\|_{X_{T}^{s,b}} \leq c_{0} \|v_{0}\|_{H^{s}(M)} +  c_{1}T^{1 - b - b'} \|u\|_{(X_{1/\sigma}^{s,b})_{T}} \|v\|_{X_{T}^{1, b}},\\
\|\Phi_{2}(u,v)\|_{(X_{1/\sigma}^{s, b})_{T}} \leq  c_{0} \|u_{0}\|_{H^{s}(M)} +  c_{1}T^{1 - b- b'} \| v \|_{X_{T}^{s, b}} \| v \|_{X_{T}^{1, b}}. \\
\end{cases}
\end{equation}

Consider $\Phi(u,v):=(\Phi_{1}(u,v), \Phi_{2}(u,v))$. Then, from   \eqref{BCL3}, we get
\begin{equation}\label{MP7}
\|\Phi(u,v)\|_{ \mathcal{X}^{s,b}_{T}} \leq c_{0} \|(v_0,u_0)\|_{\bH^{s}(M)} + c_{1} T^{1- b - b'} \|(v,u)\|_{\mathcal{X}_{T}^{s,b}} \|(v, u)\|_{\mathcal{X}_{T}^{1, b}}.
\end{equation}

If  $(v,u) \in B^{R}_{T}$ we must choose  $T$ so that 
\[
c_{1} T^{1- b - b'} \|(v,u)\|_{\mathcal{X}_{T}^{1,b}} < \frac{1}{2}.
\]
that is,
\[
T =  \left(\frac{1}{ 4 c_{1}  \|(v,u)\|_{\mathcal{X}_{T}^{1,b}}}\right)^{\frac{1}{1 - b - b'}} > 0.
\]
By the continuous dependence of $H^1(M)$, it follows that 
\[
 \|(v,u)\|_{\mathcal{X}_{T}^{1,b}} \leq C \|(v_0,u_0)\|_{\bH^{1}}.
\]
Thus, we have
\[
\left( \frac{1}{4 c_1 C \|(v_0,u_0)\|_{\bH^{1}} }\right)^{\frac{1}{1 - b - b'}} \leq  \left(\frac{1}{ 4 c_{1}  \|(v,u)\|_{\mathcal{X}_{T}^{1,b}}}\right)^{\frac{1}{1 - b - b'}} = T,
\]
 such that the solution  exists on $[0, T]$,  and
\[
 \|(v,u)\|_{L^{\infty}((0, T); \bH^{s})} \lesssim  \|(v,u)\|_{\mathcal{X}_{T}^{s,b}} \lesssim \|(v_0,u_0)\|_{\bH^{s}}.
\]

Observe that  $T$ depends only on the  $\bH^1$ norm of the initial data. Moreover, by
 \eqref{MP5}, there exists a fixed constant   $C = C(\sigma, A_{opt}, B, \alpha, d) > 0$ such that 
\[
0 < C(\sigma, A_{opt}, B, \alpha, d) \leq  \frac{1}{\|(v(t), u(t))\|_{\bH^{1}}}.
\]
 Thus, $T$ is bounded from below by a positive constant, and we can iterate this argument to extend  the solution 
$(v, u)$ on any time interval  $[0, T]$   for   $s > 1$. 
\end{dem}
\subsection*{Acknowledgment}
The authors would like to thank 
the anonymous referee for the careful reading of the manuscript and
many constructive comments and suggestions that considerably improved the
presentation.

\bibliographystyle{plain}

\end{document}